\documentclass[a4paper,11pt]{article}
\usepackage{geometry,amssymb,amsmath,enumerate,latexsym,array}
\usepackage[dvips]{graphicx}
\usepackage[all]{xy}
\newcommand{\id}{i\!d}
\newcommand{\ep}{\varepsilon}
\newcommand{\bA}{\mathsf{A}}
\newcommand{\bB}{\mathsf{B}}
\newcommand{\bC}{\mathsf{C}}

\newcommand{\bN}{\mathbb{N}}
\newcommand{\bP}{\mathsf{P}}
\newcommand{\sets}{\mathsf{Sets}}

\newcommand{\lra}{\leftrightarrow}
\newcommand{\lan}{\langle}
\newcommand{\ran}{\rangle}

\def\ob{{\mathrm{ Ob}}}

\def\bar{\overline}

\def\leq{\leqslant}
\def\geq{\geqslant}
\newcommand{\arr}{\mathrm{Arr}}
\newtheorem{example}{Example}[section]
\newtheorem{Def}[example]{Definition}
\newtheorem{rem}[example]{Remark}
\newtheorem{prop}[example]{Proposition}
\newtheorem{cor}[example]{Corollary}
\newtheorem{alg}[example]{Algorithm}
\newtheorem{thm}[example]{Theorem}
\newtheorem{lem}[example]{Lemma}

\newenvironment{proof}{\noindent {\bf Proof} }{ \hfill
$\Box$ \mbox{}}
\begin{document}

\title { Using rewriting systems to compute Kan extensions\\
and induced actions of categories\thanks{
AMS 1991 CLASSIFICATION: 68Q42,18A40,68Q40,18-04}}
\author{ Ronald Brown\thanks{Research partially supported by INTAS
Project 94-436 ext `Algebraic K-theory, groups and categories'.}
\\School of Mathematics \\ University of Wales, Bangor \\
Gwynedd LL57 1UT \\ United Kingdom \\r.brown@bangor.ac.uk
\and Anne Heyworth\thanks{ Supported 1995-8 by an EPSRC
Earmarked Research Studentship, `Identities among relations for
monoids and categories', and 1998-9 by  a University of Wales, Bangor,
Research Assistantship.}\\ School of Mathematics \\ University of Wales, Bangor \\
Gwynedd LL57 1UT \\ United Kingdom \\ map130@bangor.ac.uk
 } \maketitle

\begin{abstract}
The basic method of rewriting for words in a free monoid given a
monoid presentation is extended to rewriting for paths in a free
category given a `Kan extension presentation'. This is related to
work of Carmody-Walters on the Todd-Coxeter procedure for Kan
extensions, but allows for the output data to be infinite,
described by a language. The result also allows rewrite methods to
be applied in a greater range of situations and examples, in terms
of induced actions of monoids, categories, groups or groupoids.
\end{abstract}

\newpage
\section{Introduction}

This paper extends the usual  rewriting procedures for words $w$
in a free monoid to terms $x|w$ where $x$ is an element of a set
and $w$ is a word. Two kinds of rewriting are involved here. The
first is the familiar $x|ulv \to x|urv$ given by a relation
$(l,r)$. The second derives from a given action of certain words
on elements, so allowing rewriting $x|F(a) v \to x \cdot a|v$ (a
kind of tensor product rule). Further, the elements $x$ and $x
\cdot a$ are allowed to belong to different sets.

The natural setting for this rewriting is a {\em presentation}\ of
the form  $kan \lan \Gamma | \Delta | RelB | X | F \ran$ where:
\begin{itemize}
  \item $\Gamma,\Delta $ are (directed) graphs;
  \item $X: \Gamma \to \sets$ and $F: \Gamma \to P \Delta$ are graph
morphisms to the category of sets and the free category on
$\Delta$ respectively;
  \item  and $RelB$ is a set of relations on the free category $P\Delta$.
\end{itemize}
 The main result defines rewriting procedures on the $P
\Delta $-set $$T:= \bigsqcup_{B \in \ob\Delta} \bigsqcup_{A \in
\ob \Gamma} XA \times P\Delta(FA,B). $$ When such rewriting
procedures complete, the associated normal form gives in effect a
computation of what we call the {\em Kan extension defined by the
presentation}.

So the power of rewriting theory may now be brought to bear on a
much wider range of combinatorial enumeration problems.
Traditionally rewriting is used for solving the word problem for
monoids. It has also been used for coset enumeration problems
\cite{Redfern,Hurt}.  It may now also be used in the specification of
\begin{enumerate}[i)]
\item equivalence classes and equivariant equivalence classes,
\item arrows of a category or groupoid,
\item right congruence classes given by a relation on a monoid,
\item orbits of an  action of a group or monoid.
\item conjugacy classes of a group,
\item coequalisers, pushouts and colimits of sets,
\item induced permutation representations of a group or monoid.
\end{enumerate}
and many others.

In this paper we are concerned with the description of the theory
and the implementation in {\sf GAP} of the procedure with respect
to one ordering. It is hoped to consider implementation of
efficiency strategies and other orderings on another occasion. The
advantages of our abstraction should then become even clearer,
since one efficient implementation will be able to apply to a
variety of situations, including some not yet apparent.

We would like to acknowledge the help given by  Larry Lambe in
computational and mathematical advice  since the early 1990s. He
further suggested in 1995 that data structures of free categories
implemented by Brown and Dreckmann could be relevant to work of
Carmody and Walters on computations of Kan extensions. In visits
in 1996 and 1997 under an EPSRC Visiting
Fellowship\footnote{`Symbolic computation and Kan extensions',
GR/L22416, 1996-7.} he gave further crucial direction to the work,
including suggestions on the connections with Gr\"obner bases
which are developed elsewhere

The papers \cite{BLW,CaWa1, CaWa2, Rosebrugh} were very influential
on the current work.
\section{Kan Extensions of Actions}

The concept of the Kan extension of an action  will be defined in
this section with some familiar examples to motivate the
construction listed afterwards.

There are two types of Kan extension (the details are in Chapter
10 of \cite{Mac}) known as right and left. Which type is right and
which left varies according to authors' chosen conventions. In
this text only one type is used (left according to \cite{CaWa1},
right according to other authors) and to save conflict it will be
referred to simply as ``the Kan extension'' - it is the colimit
one, so there is an argument for calling it a co-Kan, and the
other one simply Kan, but we shall not presume to do that here.\\

Let $\bA$ be a category. A \textbf{category action} $X$ of $\bA$
is a contravariant functor $X:\bA \to \sets$. This means that for
every object $A$ there is a set $XA$ and the arrows of $\bA$ act
on the elements of the sets associated to their sources to return
elements of the sets associated to their targets. So if $a_1$ is
an arrow in $\bA(A_1,A_2)$  then $XA_1$ and $XA_2$ are sets and
$Xa_1 : XA_1 \to XA_2$ is a function where $Xa_1(x)$ is denoted $x
\cdot a_1$. Furthermore, if $a_2\in\bA(A_2,A_3)$ is another arrow
then $(x\cdot a_1)\cdot a_2=x.(a_1a_2)$ so the action preserves
the composition. This is equivalent to the fact that
$Xa_2(Xa_1(x))=X(a_1a_2)(x)$ i.e. $X$ is a contravariant functor.
The action of identity arrows is trivial, so if $\id$ is an
identity arrow at $A$ then $x \cdot \id = x$ for all $x \in XA$.

Given the category $\bA$ and the action defined by $X$, let $\bB$
be a second category and let $F:\bA \to \bB$ be a covariant
functor. Then an \textbf{extension of the action $X$ along $F$} is
a pair $(K,\ep)$ where $K:\bB \to \sets$ is a contravariant
functor and $\ep:X \to F \circ K$ is a natural transformation.
This means that $K$ is a category action of $\bB$ and $\ep$ makes
sure that the action defined is an extension with respect to $F$
of the action already defined on $A$. So $\ep$ is a collection of
functions, one for each object of $\bA$, such that
$\ep_{src(a)}(Xa)$ and $K(F(a))$ have the same action on elements
of $K(F(src(a))$.\\

The \textbf{Kan extension of the action $X$ along $F$} is an
extension $(K,\ep)$ of the action  with the universal property that
for any other extension of the action $(K',\ep')$ there exists a
unique natural transformation $\alpha:K \to K'$ such that
$\ep'=\alpha \circ \ep$. Thus $K$ may thought of as the universal
extension of the action of $\bA$ to an action of $\bB$.
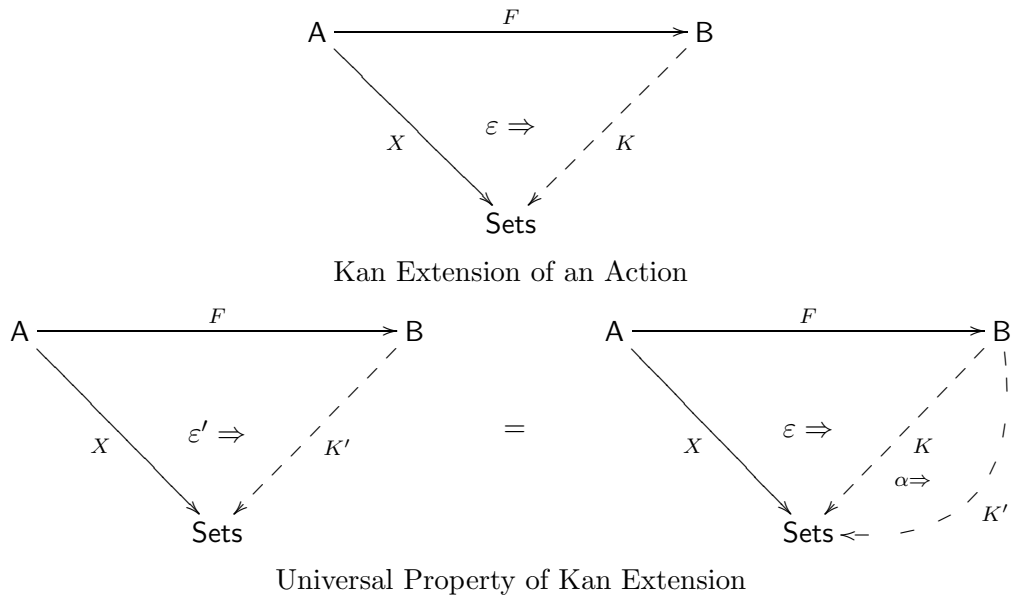
$$
\xymatrix{ {\bA} \ar[rrrr]^F \ar[ddrr]_X &&&& {\bB} \ar@{-->}[ddll]^K \\
           && \ep \Rightarrow && \\
           && {\sets} && \\ }
$$$$\text{Kan Extension of an Action}
$$
$$
\xymatrix{ {\bA} \ar[rrrr]^F \ar[ddrr]_X &&&& {\bB} \ar@{-->}[ddll]^{K'}
           &&  {\bA} \ar[rrrr]^F \ar[ddrr]_X &&&& {\bB} \ar@{-->}[ddll]^K
             \ar@/^3pc/@{-->}[ddll]^{K'} \\
           && \ep' \Rightarrow &&
           & = &  && \ep \Rightarrow &  \ar@{}[d]|{\alpha\Rightarrow}&\\
           && {\sets} &&  &&  && {\sets} && \\ }
$$$$
\text{Universal Property of Kan Extension}
$$
\section{Examples}
Mac\,Lane wrote in section 10.7 of \cite{Mac} (entitled ``All
Concepts are Kan Extensions'') that ``the notion of Kan extensions
subsumes all the other fundamental concepts of category theory'' .
We now illustrate his statement by showing how some familiar
problems can be expressed in these terms, and will later see how
our computational methods apply to these problems. Most of these
examples are also familiar from  \cite{CaWa1,Rosebrugh}.
Throughout these examples we use the same notation as the
definition, so the pair $(K,\ep)$ is the Kan extension of the
action $X$ of $\bA$ along the functor $F$ to $\bB$. By a monoid
(or group) ``considered as a category'' we mean the one object
category with arrows corresponding to the monoid elements and
composition defined by composition in the monoid.

\textbf{1) Groups and Monoids}\\ Let $\bB$ be a monoid regarded as
a category with one object $0$. Let $\bA$ be the singleton category,
acting trivially on a one point set $X0$, and let $F:\bA \to
\bB$ be the inclusion map. Then the set $K0$ is isomorphic
to the set of elements of the monoid and the right action of the
arrows of $\bB$ is right multiplication by the monoid elements.
The natural transformation maps the unique element of $X0$
to the element of $K0$ representing the monoid identity.

\textbf{2) Groupoids and Categories}\\ Let $\bB$ be a category.
Let $\bA$ be the (discrete) sub-category of objects of $\bB$ with
identity arrows only. Let $X$ define the trivial action of $\bA$
on a collection of one point sets $\bigsqcup_B  XB$ (one for each
object $B$ of $\bB$), and let $F:\bA\to\bB$ be the inclusion map.
Then the set $KB$ for $B\in\ob\bB$ is isomorphic to the set of
arrows of $\bB$ with target $B$ and the right action of the arrows
of $\bB$ is defined by right composition. The natural
transformation $\ep$ maps the unique element of a set $XB$ to the
representative identity arrow for the object $FB$ for every
$B\in\ob\bA$.\\

\textbf{3) Cosets, and Congruences on Monoids}\\ Let $\bB$ be a
group considered as a category with one object $0$, and let
$F:\bA\to \bB$ be the inclusion of the  subgroup $\bA$. Let $X$
map the object of $\bA$ to a one point set. The set $K0$
represents the (right) cosets of $\bA$ in $\bB$, with the right
action of any group element  $b$ of $\bB$  taking the
representative of the coset $\bA g$ to the representative of the
coset $\bA gb$. The left cosets can be similarly represented,
defining the right action $K$ by a left action on the cosets. The
natural transformation $\ep$ picks out the representative for the
subgroup $\bA$.
\\
Alternatively, let $\bB$ be a monoid considered as a category with
one object $0$ and let $\bA$ be generated by arrows which map under
$F$ to a set of generators for a right congruence. Then the set
$K0$ represents the congruence classes, the action of any arrow
$b$ of $\bB$ (monoid elements) taking the representative (in $K0$)
of the class $[m]$ to the representative of the class $[mb]$. The
natural transformation picks out the representative for the class
$[id]$. (As above, left congruence classes may also be expressed
in terms of a Kan extension.)\\

\textbf{4) Orbits of Group Actions}\\ Let $\bA$ be a group thought
of as a category with one object $0$ and let $X$ define the action
of the group on a set $X0$. Let $\bB$ be the trivial category on
the object $0$ and let $F$ be the null functor. Then the set $K0$
is a set of representatives of the distinct orbits of the action
of $\bA$ and the action of $\bB$ on $K0$ is trivial. The natural
transformation $\ep$ maps each element of the set $X0$ to its
orbit representative in $K0$.\\

\textbf{5) Colimits in $\sets$}\\ Let $X:\bA\to \sets$ be any
functor on the small category $\bA$ and let $F:
\bA\to\bB$ be the null functor to the trivial category. Then the
Kan extension corresponds to the colimit of (the diagram)
$X:\bA\to\sets$; $K0$ is the colimit object, and $\ep$ defines the
colimit functions from each set $XA$ to $K0$. Examples of this
are: (i) when $\bA$ has two objects $A_1$ and $A_2$, and two
non-identity arrows $a_1,a_2 :A_1 \to A_2$;  the colimit is then
the \emph{coequaliser} of the functions $Xa_1$ and $Xa_2$ in
$\sets$; (ii) when $\bA$ has three objects $A_1$, $A_2$ and $A_3$
and two arrows $a_1:A_1 \to A_2$ and $a_2:A_1 \to A_3$; the
colimit is then the \emph{pushout} of the functions $Xa_1$ and
$Xa_2$ in $\sets$.\\

\textbf{6) Induced Permutation Representations}\\ Let $F:\bA \to
\bB$ be a morphism of groups, thought of as a functor of
categories. Let  $X$ be a right action of the group $\bA$ on the
set $X0$. The Kan extension of the action along $F$ is known as
the action of $\bB$ \emph{induced} from that of $\bA$ by $F$; it
is sometimes written $F_*(X)$. There are simple methods of
constructing the set $K 0$ in this case. For example if $F$ is
surjective, then $F_*X$ may be taken to be the set $X$ factored by
the action of $ker( F)$, while if $F$ is injective then  $F_*X$
may be taken to be the set $X\times S$ where $S$ is a transversal
of $F(\bA)$ in $\bB$, with an appropriate action ref{appropact}.
A corresponding
description of the Kan extension is more difficult for monoid
actions.\\

This last example is very close to the full definition of a Kan
extension. A Kan extension \emph{is} the action of the category
$\bB$ induced from the action of $\bA$ by $F$ together with $\ep$
which shows how to get from the $\bA$-action to the $\bB$-action.
The point of giving the other examples is to show that Kan
extensions can be used as a method of representing a variety of
situations.

\section{Presentations of Kan Extensions of Actions}
The problem that has been introduced is that of ``computing a Kan
extension''. In order to keep the analogy with computation and
rewriting for presentations of monoids we propose a definition of
a {\em presentation} of a Kan extension.

First, we set out our notation for free categories. Let $\Delta$
be a directed graph, that is $\Delta$ consists of two functions
$src,tgt: \arr\Delta \to \ob\Delta$. Any small category $\bP$ has
an underlying graph $U\bP$.  The {\em free category} $P\Delta$ on
$\Delta$ consists of the objects of $\Delta$ with an identity arrow at
each object and non identity arrows
$p:B \to B'$ given by the sequences $(d_1,d_2, \ldots, d_n)$ of arrows of
$\Delta$ which are composable, i.e. $tgt (d_i)= src(d_{i+1}),
1=1,\ldots,n-1$, and such that $src(d_1)=B, tgt(d_n)=B'$.
As
usual, such a word is written $d_1\ldots d_n:B \to B'$, and
composition is by juxtaposition. Of course the free functor $P$ is
left adjoint to the forgetful functor $U$.

A {\em graph  of relations } $Rel$ for the free category $P\Delta$
has objects those of $\Delta$ and arrows $B \to B'$ a set of pairs
$(l,r)$ such that $l,r:B \to B'$ in $\Delta$. Then the quotient
category $P\Delta/Rel$ is defined.

A {\em presentation } $cat\lan \Delta|Rel\ran$ for a category
$\bB$ consists of a graph $\Delta$ of generators of $\bB$ and a
graph of relations for $P\Delta$  such that the natural morphism
of categories $P\Delta \to \bB$ induces an isomorphism of
categories $(P\Delta)/Rel \to \bB$.
(For an introduction to category presentations see \cite{Mit}).

Next, we define `Kan extension data'.

\begin{Def}
A \textbf{Kan extension data} $(X',F')$ consists of small
categories $\bA$, $\bB$ and functors $X':\bA\to\sets$ and
$F':\bA\to\bB$.
\end{Def}
\begin{Def}
  A \textbf{Kan extension presentation} is a quintuple
$\mathcal{P}:=kan\lan \Gamma|\Delta|RelB|X|F \ran$ where
\begin{enumerate}[i)]
\item
$\Gamma$ and $\Delta$ are graphs,
\item
$cat\lan \Delta | RelB \ran$ is a category presentation,
\item
$X: \Gamma \to U \sets$ is a graph morphism,
\item
$F: \Gamma \to U P\Delta$ is a graph morphism.
\end{enumerate}

We say $\mathcal{P}$ {\em presents} the Kan extension data
$(X',F')$ where $X':\bA\to\sets$ and $F':\bA\to\bB$ if
\begin{enumerate}[i)]
\item
$\Gamma$ is a generating graph for $\bA$ and $X:\Gamma\to\sets$ is the
restriction of $X':\bA\to\sets$,
\item
$cat\lan \Delta|RelB\ran$ is a category presentation of $\bB$,
\item
$F:\Gamma\to P\Delta$ induces $F':\bA\to\bB$.
\end{enumerate}
We also say $\mathcal{P}$ \textbf{presents} the Kan extension
$(K,\ep)$ of the Kan extension data $(X',F')$. The presentation is
\textbf{finite} if $\Gamma$, $\Delta$ and $RelB$ are finite.
\end{Def}

\begin{rem}\emph{
The fact that $X,\,F$ induce $X',\,F'$ implies extra conditions on
$X,\,F$ in relation to $\bA$ and $\bB$. In practice we need only
the values of $X',\,F'$ on $\Gamma$.  In other words, a given  Kan
extension presentation always defines a Kan extension data where
$\bA$ is the free category $P\Gamma$ and $(X',F')$ are induced by
$X,F$. This is analogous to the fact that for coset enumeration of a
subgroup $H$ of $G$ where $G$ has presentation $grp\lan
\Delta|R\ran$ we need only that $H$ is generated by certain words
in the set $\Delta$. }\end{rem}

\section{$\bP$-sets}

In this section we establish the concepts and notation used
to apply rewriting procedures to presentations
of Kan extensions of actions. Our terminology is modelled on that
standard in rewriting theory.

\begin{Def}
\label{Pset} Let  $\bP$ be a category. A \textbf{$\bP$-set} is a
set $T$ together with a function $\tau:T \to \ob\bP$ and a partial
action $\cdot$ of the arrows of $\bP$ on $T$. The action satisfies
the following properties for all $t\in T, p,q \in \arr \bP$:
\begin{enumerate}[i)]
  \item if $\tau(t)=src(p) $ then $t\cdot p$ is defined and
  $\tau(t\cdot p) = tgt(p)$;
  \item $t \cdot \id_{\tau(t)} = t$;
  \item  $(t\cdot p)\cdot q=t \cdot (pq)$ if the left hand side is
  defined.
\end{enumerate}
\end{Def}

\begin{Def}
\label{redrel} A \textbf{reduction relation} on a $\bP$-set $T$ is
a relation $\to$ on $T$ such that for all $t_1,t_2 \in T$, $t_1
\to t_2$ implies $\tau(t_1)=\tau(t_2)$. The reduction relation
$\to$ on the $P$-set $T$ is \textbf{admissible} if for all $t_1,
t_2 \in T$, $t_1 \to t_2$ implies $t_1 \cdot p \to t_2 \cdot p$
for all $p \in \arr\bP$ such that $src(p)=\tau(t_1)$.
\end{Def}

For the rest of this paper we assume that $\mathcal{P}=kan \lan
\Gamma | \Delta | RelB |X |F \ran$ is a presentation of a Kan
extension. The following definitions will be used throughout. Let
$\bP$ denote the free category $P\Delta$. Then define
\begin{equation} \label{Tdef}
T := \bigsqcup_{B\in\ob\Delta}\bigsqcup_{A\in\ob\Gamma} XA \times
\bP(FA, B)
\end{equation}
The elements of the set $T$ will be referred to as {\em terms},
and a pair  $(x,p) \in XA \times \bP(FA,B)$ will be written $x|p\,$.
 The function $\tau:T\to\ob\bP$ is defined by
$$\tau(x|p):=tgt(p) \text{ for } x|p \in T.$$ Then $T$ becomes a
$\bP$-set by the action  $$(x|p) \cdot q := x|pq \text{ for } x|p
\in T, \ q \in \arr\bP \text{ when } src(q)=\tau(x|p).$$

A \textbf{rewrite system} for a Kan presentation $\mathcal{P}$ is
a pair $R:=(R_T,R_P)$  such that
\begin{enumerate}[(i)]
  \item $R_T$ is a reduction relation on the $\bP$-set $T$;
  \item  $R_P$ is a set of relations on $\bP$, so that $(l,r)
  \in R_P$ implies $l,r\in \bP(B,B')$ for some $B,B' \in \ob(\Delta)$.
\end{enumerate}

The \textbf{initial rewrite system} that results from the
presentation is the pair $R_{init}:=(R_\ep,  R_K)$ defined as follows.
\begin{align*}
R_\ep:&=\{ (x|Fa,x\cdot a|id_{FA_2}) | x\in XA_1, a \in
\Gamma(A_1, A_2), A_1, A_2 \in \ob \Gamma \}. \\ R_K:&=RelB.
\end{align*}

The first type of rule we call the `$\ep$-rules' $R_\ep \subseteq
T \times T$. They are to ensure that the action is an extension by
$F$ of the action of $P\Gamma$
 -- this is the requirement for $\ep:X \to KF$ to be a natural
transformation.

The second  type we call the `$K$-rules' $R_K \subseteq
\arr \bP \times \arr \bP.$ They are to ensure that the action
preserves the relations and so gives a functor on the quotient
$\bB=(P\Delta)/RelB$.
\begin{rem}\emph{
If the Kan extension presentation is finite then $R_{init}$ is finite.
The number of initial rules is by definition
$(\Sigma_{a\in\arr\Gamma}|X(src(a))|)+|RelB|$. }
\end{rem}

\begin{Def}
The \textbf{reduction relation $\to_R$ generated by} a rewrite
system $R=(R_T, R_P)$ on the $\bP$-set $T$ is defined as $t_1
\to_R t_2$ if and only if one of the following is true:
\begin{enumerate}[i)]
\item
There exist $(s,u)\in R_T, q \in \arr \bP$
such that $t_1= s \cdot q$ and $t_2= u \cdot q$.
\item
There exist $(l,r) \in R_P$, $s \in T$, $q \in \arr \bP$ such that $t_1
=s \cdot l q$ and $t_2= s\cdot rq$.
\end{enumerate}
Then we say $t_1$ {\em reduces to} $t_2$ by the rule $(s,u)$ or by
$(l,r)$ respectively.
\end{Def}

Note that $\to_R$ is an admissible reduction relation on $T$.
The relation
$\stackrel{*}{\to}_R$ is defined to be the reflexive, transitive
closure of $\to_R$ on $T$, and $\stackrel{*}{\lra}_R$ is the
reflexive, symmetric, transitive closure of $\to_R$.

\begin{rem}\emph{
Essentially, the rules of $R_P$ are two-sided and apply to any
substring to the right of the separator $|$. This distinguishes
them from the one-sided rules of $R_T$ -- these might be called
`tagged rewrite rules', the `tag' being the part $x$ to the left
of the separator of $x|p$, but in a more general sense than
previous uses since the tags are being rewritten.}
\end{rem}

\begin{lem}
Let $R$ be a rewrite system on a $\bP$-set $T$. Then
$\stackrel{*}{\lra}_R$ is an admissible equivalence relation on the
$\bP$-set $T$.
\end{lem}

The proof is straightforward.

The equivalence class of $t \in T$ under $\stackrel{*}{\lra}_R$
will be denoted $[t]$. A suggestive notation for the class $[x|p]$
would also be $x \otimes p$. \\

We apply the standard terminology of reduction relations
to the reduction relation $\to_R$ on $T$. In particular we have a
notion of $\to_R$ being complete. A rewrite system $R :=
(R_T, R_P)$ will be called  \textbf{complete} when $\to_R$ is complete. In
this case $\stackrel{*}{\lra}_R$ admits a normal form function.

We expect that a Kan extension $(K, \ep)$ is given by a set $KB$
for each $B \in \ob \Delta$ and a function $Kb: KB_1 \to KB_2$ for
each $b:B_1 \to B_2 \in \bB$  (defining the functor $K$) together
with a function $\ep_A: XA \to KFA$ for each $A \in \ob \bA$ (the
natural transformation). This information can be given in four
parts:
\begin{itemize}
  \item  the set $\bigsqcup_B  KB$;
  \item  a function $\bar{\tau}:\bigsqcup_B KB \to \ob \bB$;
  \item  a partial function (action) $\bigsqcup_B KB \times \arr\bP \to \bigsqcup_B
  KB$;
  \item  and a function $\ep: \bigsqcup_A XA \to \bigsqcup_B  KB$.

\end{itemize}
Here $\bigsqcup_B  KB$ and $\bigsqcup_A XA$ are the disjoint
unions of the sets $KB$, $XA$ over $\ob \bB$, $\ob \bA$
respectively; if $z \in KB$ then $\bar{\tau}(z)=B$
 and if further $src(p)=B$ for $p \in \arr\bP$ then $z \cdot p$
is defined.

\begin{thm}
Let $\mathcal{P}=kan \lan \Gamma | \Delta | RelB | X F \ran$ be a
Kan extension presentation, and let $\bP$, $T$, $R=(R_\ep, R_K)$
be defined as above. Then the Kan extension $(K,\ep)$ presented by
$\mathcal{P}$ may be  given by the following data:
\begin{enumerate}[i)]
\item
the set $\bigsqcup_B KB = T/ \stackrel{*}{\lra}_R$,
\item
the function $\bar{\tau}: \bigsqcup_B KB \to \ob \bB$ induced by
$\tau:T \to \ob \bP$,
\item
the action of $\bB$ on $ \bigsqcup_B KB$ induced by the action of
$\bP$ on $T$,
\item
the natural transformation $\ep$ determined by $x \mapsto [x|
\id_{FA}]$ for $x \in XA, \;A \in \ob \bA$.
\end{enumerate}
\end{thm}

\begin{proof}
We give the proof in some detail since this is helpful  for the
implementations described in the next section.

\textbf{Claim}
$\stackrel{*}{\lra}$ preserves the function $\tau$.

\begin{proof}
 We prove that $\lra$, the
symmetric closure of $\to$, preserves $\tau$. Let $t_1,t_2\in T$
so that  $t_1 \lra t_2$. From the definition of $\to$ there are
two possible situations. For the first case suppose that there
exist $(s_1,s_2) \in R_\ep$ such that $t_1 = s_1 \cdot q$ and $t_2
= s_2 \cdot q$ for some $q \in \arr \bP$. Clearly
$\tau(t_1)=\tau(t_2)$. For the other case suppose that there exist
$(l,r)\in R_K$ such that $t_1=s \cdot(lq)$ and $t_2 = s\cdot(rq)$
for some $s \in T$, $q \in \arr \bP$. Again, it is clear that
$\tau(t_1)=\tau(t_2)$. Hence $\bar{\tau}:T/\stackrel{*}{\lra}_R \;
\to \ob\bP$ is well-defined by $\bar{\tau}[t]=\tau(t)$.
\mbox{ } \end{proof}

\textbf{Claim}
$T/\stackrel{*}{\lra}$ is a $\bB$-set.

\begin{proof}
First we prove that $\bB$ acts on the equivalence classes of $T$
with respect to $\stackrel{*}{\lra}$. An arrow of $\bB$ is an
equivalence class $[p]$ of arrows of $\bP$ with respect to $RelB$.
It is required to prove that $[t]\cdot p:=[t\cdot p]$ is a well
defined action of $\bP$ on $T/\stackrel{*}{\lra}$ such that
$[t]\cdot p=[t]\cdot q$ for all $p=_{RelB}q$. Let $t\in T,
p\in\arr\bP$ be such that $\tau[t]=src[p]$ i.e. $\tau(t)=src(p)$.
Then $t\cdot p$ is defined. Suppose $s \stackrel{*}{\lra} t$. Then
$[s\cdot p]=[t\cdot p]$ since $s\cdot p \stackrel{*}{\lra} t\cdot
p$, whenever $s\cdot p, t\cdot p$ are defined. Suppose $p =_{RelB}
q$. Then $[t\cdot p]=[t\cdot q]$ since $t\cdot p\,
{\stackrel{*}{\lra}}_{R_K}\, t\cdot q$, whenever $t\cdot p, t\cdot
q$ are defined and $({\stackrel{*}{\lra}}_{RelB})\, \subseteq\,
(\stackrel{*}{\lra})$. Therefore $\bP$ acts on
$T/\stackrel{*}{\lra}$. This action preserves the relations of
$\bB$ and so defines an action of $\bB$ on $T/\stackrel{*}{\lra}$.
Furthermore $\bar{ \tau }( [ t ] \cdot p ) = \bar{\tau}[t \cdot p
] = tgt( p )$ and if $q \in \bP$ such that $src(q)=tgt(p)$ then $(
[ t ] \cdot p ) \cdot q = [ ( t \cdot p ) \cdot q ] = [ t \cdot (
pq ) ] = [ t ] \cdot pq$.
\end{proof}

  The Kan extension may now be defined. For $B\in\ob\bB$ define
\begin{equation}
\label{Kobjects}
KB:=\{ [x|p] : \bar{\tau}[x|p] = B \}.
\end{equation}
For $b:B_1\to B_2$ in $\bB$ define
\begin{equation}
\label{Karrows}
Kb:KB_1\to KB_2 : [t] \mapsto [t\cdot p] \text{ for }[t]\in KB_1\text{
where } p \in [b].
\end{equation}
It can be verified that this definition of the action is a functor
$K:\bB \to \sets$. Then define
\begin{equation}
\label{nattran}
\ep:X\to KF:x\mapsto[x|\id_{FA}]\text{ for }x\in XA,A\in\ob\bA.
\end{equation}
It is straightforward to verify that this is a natural
transformation. Therefore $(K,\ep)$ is an extension of the action
$X$ of $\bA$. The proof of the universal property of the extension
is as follows. Let $K':\bB\to\sets$ be a functor and $\ep':X \to
K'F$ be a natural transformation. Then  $\alpha:K \to K'$, defined
by $$\alpha_B[x|p] = K'(f) (\ep_A'(x))  \text{ for } [x|p]\in
KB,$$ is a  natural transformation which  satisfies $\ep \circ \alpha = \ep'$
 and is clearly the only such.
\end{proof}

\section{Rewriting Procedures for Kan Extensions}

In this section we will explain the completion process for the
initial rewrite system. To this end we give a convenient
 notation for the implementation of the data structure for a
\emph{finite} presentation $\mathcal{P}$ of a Kan extension.
The functions which work with this structure form a package {\sf Kan}
which is being submitted as a share package for {\sf GAP}.
\subsection{Input Data}
In the {\sf GAP} system, a symbol such as $b_3$ can be defined only as a
`generator'. This explains the use of the term `generator' in the
following.

\begin{enumerate}[]
\item
$\mathtt{ObA}$ \quad This is a list $[1,2,\ldots]$ of $|\ob
\Gamma|$ integers $i$ such that $i$ labels the object $A_i$ of
$\Gamma$.
\item
$\mathtt{ArrA}$ \quad This is a list of pairs of integers
$[[i_1,j_1],[i_2,j_2],\ldots]$, one for each arrow $a_k : A_{i_k}
\to A_{j_k}$ of $\arr\Gamma$. The first element of each pair is
the source of the arrow it represents, and the other entry is the
target.
\item
$\mathtt{ObB}$ \quad Similarly to $\ob\Gamma$, this is a list of
integers representing the objects of $\Delta$.
\item
$\mathtt{ArrB}$ \quad This is a list of triples
$[[b_1,i_1,j_1],[b_2,i_2,j_2],\ldots]$, one triple for each arrow
$b_k : B_{i_k} \to B_{j_k}$ of $\mathtt{\arr\Delta}$. The first
entry of each triple is a label for the arrow (in {\sf GAP} such a label
is a `generator'), and the other entries are integers representing
the source and target respectively. Note that the arrows of
$\Gamma$ did not have labels. The arrows of $\Delta$ will form
parts of the terms of $T$ whilst those of $\Gamma$ do not, so this
is why we have labels here and not before.
\item
$\mathtt{RelB}$ \quad This is a finite list of pairs of paths.
Each path is represented by a finite list $[b_1,b_2,\ldots,b_n]$
of labels of composable arrows of $\mathtt{\arr\Delta}$.  In {\sf GAP}
it is convenient to consider these lists as words $b_1 \ldots b_n$
in the generators that are labels for the arrows of $\Delta$.
\item
$\mathtt{FObA}$ \quad This is a list of $|\ob\Gamma|$ integers.
The $k$th entry represents the object of $\Delta$ which is the
image of the object $A_k$ under $F$.
\item
$\mathtt{FArrA}$ \quad This is a list of paths where the entry at
the $k$th position is the element of $\bP$, i.e. a path in
$\Delta$, which is the image of $a_k$ under $F$. The length of the
list is $|\arr \Gamma|$.
\item
$\mathtt{XObA}$ \quad This is a list $L$ of lists of distinct
({\sf GAP}) generators. There is one list $L[i]$  for each object $A_i$
in $\Gamma$, and $L[i]$ represents the elements of $XA_i$.
\item
$\mathtt{XArrA}$ \quad This is a list $M$ of lists of generators.
There is one list $M[k]$ for each arrow $a_k$ of $\Gamma$. It
represents the image under the action $Xa_k$ of the set
$X(src(a_k))$. Suppose $a_k : A_{i_k} \to A_{j_k}$ is the arrow at
entry $k$ in $\arr\Gamma$, and $[x_1,x_2,\ldots,x_m]$ is the $i$th
entry in $X \ob \Gamma$ (the image set $X(A_i)$). Then the $k$th
entry of $X \arr \Gamma$ is the list $[x_1\cdot a,x_2\cdot
a,\ldots,x_m\cdot a]$ where $x_i\in X(A_j)$.
\end{enumerate}
 All the above lists are finite since the Kan extension is
finitely presented.
In Section 8 we explain how to input this data.

\subsection{Lists}

 Elements of  $T$ are called {\em terms} and are
represented in the {\sf GAP} implementation by lists of  generators,
where the generators may be thought of as labels.  The first entry
in the list must be a label for an element of $XA$ for some $A \in
\ob \Gamma$. The subsequent entries will be labels for composable
arrows of $\Delta$, with the source of the first being $FA$.
Formally, an element $t \in T$ is represented by a list
$${\mathtt{List}}(x|p)=\begin{cases} [x,b_1,\ldots,b_n] & \text{if
}
p=b_1\ldots b_n, n \geq 1,\\ {[x,1_{FA}]} &\text{if } p=1_{FA}
.\end{cases} $$
 This also allows us to use list notation, so that if $
t=x|b_1\ldots b_n$ then $t[1]=x, t[i+1]=b_i, 1 \leq i \leq  n$.
Also, $\mathtt{Length}(t)$ means the number of elements in the
list corresponding to $t$ and $\mathtt{Position(ObA,A)}$ returns
the position of the element $A$ in the list $\mathtt{ObA}$.
If $t=[x|p]$ we also write
$t[2..]$ for $p$.

\pagebreak
\subsection{Initial Rules Procedure}

\begin{alg}[Initial Rules]
Given the data for a Kan presentation in the form of a record with the fields
named as above, the initial rewrite system $R_{init}:=(R_\ep,R_K)$
is determined.
\begin{enumerate}[1]
\item
(Input:) {\tt ObA, ArrA, ObB, ArrB, RelB, FObA, FArrA, XObA, XArrA}.
\item
(Procedure:) Set $R_\ep:=\emptyset$, then for each arrow $a \in \mathtt{ArrA}$,
set $\mathtt{i:=Position(ArrA,a)}$;\\
 $\mathtt{XA:=XObA[Position(ObA,a[1])]}$; $\mathtt{Xa := XArrA[i]}$;
and set $\mathtt{Fa := FArrA[i]}$. Then for each element  $\mathtt{x}$ in
$\mathtt{XA}$ , set $\mathtt{j:=Position(XA,x)}$ and
add the rule $\mathtt{[[x*Fa, Xa[j]];}$ to $R_\ep$. Set $R_K:=\mathtt{RelB}$.
\item
(Output:) $R_{init}:=R_\ep \sqcup R_K$.
\end{enumerate}
\end{alg}

\subsection{Orderings}

To work with a rewrite system $R$ on $T$ we will require certain
concepts of order on $T$. We  give properties of  orderings $>_X$
on $\bigsqcup_A XA$ and  $>_P$ on $\arr\bP$ to enable us to
construct an ordering $>_T$ on $T$ with the properties needed for
the rewriting procedures.

\begin{Def}
A binary operation $>$ on a set $S$  is called a \textbf{strict
partial ordering} if it is irreflexive, antisymmetric and
transitive.  It is called a \textbf{total ordering} if also for
all $x,y\in S$ either $x>y$ or $y>x$ or else $x=y$. An ordering
$>$ is \textbf{well-founded} on $S$  if there is no infinite
sequence $x_1>x_2>\cdots$ of elements of $S$. An ordering $>$ is a
\textbf{well-ordering}  if it is well-founded and a total
ordering.
\end{Def}

\begin{Def}
Let $>_P$ be a strict partial ordering on $\arr \bP$. It is called
a \textbf{total path ordering} if it induces a total order on
$\bP(B,B')$ for all objects $B,B'\in \bP$. It is called a
\textbf{well-ordering} if it is well-founded and a total path
ordering. The ordering $>_P$ is \textbf{admissible on $\arr\bP$}
if $$p >_Pq \ \Rightarrow \ upv >_Puqv$$ for all $u,v\in\arr\bP$
such that $upv,uqv\in\arr\bP$. An admissible well-ordering is
called a \textbf{monomial ordering}.
\end{Def}

\begin{lem}
Let $>_X$ be a well-ordering on the finite set $\bigsqcup_A XA$
and let $>_P$ be an admissible well-ordering on $\bP$. For
$t_1,t_2\in T$ define
\begin{align*}
t_1 >_T t_2 \text{ if }
 t_1[2..]& >_P t_2[2..] \text{ or } \\
 t_1[2..]& =t_2[2..]
   \text{ and }
t_1[1] >_X t_2[1].\\
\end{align*}
Then $>_T$ is an admissible well-ordering on the $\bP$-set $T$.
\end{lem}

\begin{proof}
It is straightforward to verify that irreflexivity, antisymmetry
and transitivity of $>_X$ and $>_P$ imply those properties for
$>_T$. The ordering $>_T$ is admissible on $T$ because it is made
compatible with the right action (defined by composition between
arrows on $\bP$) by the admissibility of $\>_P$ on $\arr\bP$. The
ordering is linear, since if $t_1,t_2\in T$ such that neither
$t_1>_Tt_2$ nor $t_2>_Tt_1$, it follows (by the linearity of $>_X$
and linearity of $>_P$ on $\arr\bP$) that $t_1=t_2$. That $>_T$ is
well-founded is easily verified using the fact that any infinite
sequence in terms of $>_T$ implies an infinite sequence in either
$>_X$ or $>_P$. Since $>_X$ and $>_P$ are both well-founded there
are no such sequences.
\end{proof}

The last result shows that there is scope for choosing different
orderings on $T$. The actual choice is even wider than this, and is
related to efficiency see \cite{Hurt} -- there may even be
completion with respect to one order and not another. We do not
discuss these matters here.

In this paper we work only with a `length-lexicographical ordering'
defined in the following way.

\begin{Def}[Implemented Ordering]
Let $>_X$ be any linear order on (the finite set) $\bigsqcup_A
XA$. Let $>_\Delta$ be a linear ordering on (the finite set)
$\arr\Delta$. This induces an admissible ordering $>_P$ on
$\arr\bP$ where
\begin{align*}
p>_Pq \Leftrightarrow &\;  \mathtt{Length}(p)>\mathtt{Length}(q)\\
 \text{ or }
        &\; \mathtt{Length}(p)=\mathtt{Length}(q)
 \text{ and there exists $k>0$ such that }\\
                     &\; p[i]=q[i] \text{ for all $i<k$ and $p[k]>_\Delta q[k]$}\\
\intertext{The ordering
$>_T$ is then defined as follows:}
t_1 >_T t_2 \text{ if } &\;
  \mathtt{Length}(t_1)>\mathtt{Length}(t_2)\\
\text{ or } &\; \mathtt{Length}(t_1) = \mathtt{Length}(t_2) \text{ and } t_1[1]>_Xt_2[1]\\
\text{ or } &\; \mathtt{Length}(t_1) = \mathtt{Length}(t_2) \text{ and there
exists } k \in [1..\mathtt{Length}(t_1)]\\
& \text{ such that } t_1[i]=t_2[i] \text{ for all }
i<k, \text{ and } t_1[k]>_\Delta t_2[k].
\end{align*}
\end{Def}

\begin{prop}
The definitions above give an admissible, length-non-increasing
well-order $>_T$ on the $\bP$-set $T$.
\end{prop}
\begin{proof}
It is immediate from the definition that $>_T$ is
length-non-increasing. It is straightforward to verify that $>_T$
is irreflexive, antisymmetric and transitive. It can also be seen
that $>_T$ is linear (suppose neither $t_1>_Tt_2$ nor $t_2>_Tt_1$
then $t_1=t_2$, by the definition, and linearity of $>_X$,
$>_\Delta$). It is clear from the definition that $>_T$ is
admissible on the $\bP$-set $T$ (if $t_1>_Tt_2$ then
$t_1.p>_Tt_2.p$). To prove that $>_T$ is well-founded on $T$,
suppose that $t_1>_Tt_2>_Tt_3>\ldots$ is an infinite sequence.
Then for each $i>0$ either
$\mathtt{Length}(t_i) > \mathtt{Length}(t_{i+1})$
or if
$\mathtt{Length}(t_i) = \mathtt{Length}(t_{i+1})$ and $t_i[1]>_Xt_{i+1}[1]$,
or if
$\mathtt{Length}(t_i) = \mathtt{Length}(t_{i+1})$ and there exists $k \in
[1..\mathtt{Length}(t_i)]$ such that $t_i[j]=t_{i+1}[j]$ for all $j<k$ and
$t_i[k]>_\Delta t_{i+1}[k]$. This implies that there is an
infinite sequence of type $n_1>n_2>n_3>\ldots$ of positive
integers from some finite $n_1$, or of type
$x_1>_Xx_2>_Xx_3>\ldots$ of elements of $\bigsqcup_A XA$ or else
of type $p_1 >_\Delta p_2 >_\Delta p_3 >_\Delta \ldots$ of arrows
of $\Delta$, none of which is possible as $>$, $>_X$, and
$>_\Delta$ are well-founded on $\bN$, $\bigsqcup_A XA$ and
$\arr\Delta$ respectively. Hence $>_T$ is well-founded.
\end{proof}

\begin{prop} \label{ordarrow}
Let $>_T$ be the order defined above. Then $p_1>_Pp_2 \Rightarrow s\cdot
p_1>_Ts\cdot p_2$.
\end{prop}
\begin{proof}
This follows immediately from the definition of $>_T$.
\end{proof}

\begin{rem}\emph{
The proposition can also be proved for the earlier definition of $>_T$
induced from $>_X$ and $>_P$.
}\end{rem}

\subsection{Reduction}

Now that we have defined an admissible well-ordering on $T$ it is
possible to discuss when a reduction relation generated by a
rewrite system is compatible with this ordering.

\begin{lem}
Let $R$ be a rewrite system on $T$. Orientate the rules of $R$ so that
for all $(l,r)$ in $R$, if $l,r\in\arr\bP$ then $l >_P r$ and if $l,r\in
T$ then $l >_T r$.
Then the reduction relation $\to_R$ generated by $R$ is compatible with
$>_T$.
\end{lem}
\begin{proof}
Let $t_1,t_2\in T$ such that $t_1\to_Rt_2$.
There are two cases to be considered, by Definition \ref{redrel}.
For the first case let $t_1=s_1\cdot p$, $t_2=s_2\cdot p$ for some
$s_1,s_2\in
T$, $p\in\arr\bP$ such that $(s_1,s_2)\in R$. Then $s_1>_Ts_2$. It
follows that
$t_1>_Tt_2$ since $>_T$ is admissible on $T$.
For the second case let $t_1=s \cdot p_1q$, $t_2=s \cdot p_2q$ for some
$s \in T$, $p_1,p_2,q\in\arr\bP$ such that $(p_1,p_2)\in T$. Then
$p_1>_Pp_2$
and so by Proposition \ref{ordarrow} $s\cdot p_1>_Ts\cdot p_2$. Hence
$t_1>_Tt_2$ by
admissibility of $>_T$ on $T$. Therefore, in either case $t_1>_Tt_2$ so
$\to_R$ is compatible with $>_T$.
\end{proof}

It is a standard result that if a reduction relation is compatible
with an admissible well-ordering, then it is Noetherian. The next
algorithm describes the function $\mathtt{Reduce}$.

\begin{alg}[Reduce]
Given a term $t \in T$ and a rewrite system $R =(R_P,R_P)$ a term
$t_n \in [t]$, which is irreducible with respect to $\to_R$, is determined.
\begin{enumerate}
\item
(Input:) A term $t$ (as a list) and a rewrite system $R$ (as a list of pairs of
lists).
\item
(Loop:) While any left hand side of any pair occurs as a sublist of $t$ replace
that part of $t$ with the right hand side to define a reduced term $t'$.
Repeat until no left hand side of any pair
occurs in the reduced term $t'$.
\item
(Output:) A term $t'$ that is irreducible with respect to $\to_R$.
\end{enumerate}
\end{alg}

\subsection{Critical Pairs}

We can now discuss what properties of $R$ will make $\to_R$ a complete
(Noetherian and confluent) reduction relation. By standard abuse of
notation the rewrite system $R$ will be called complete when $\to_R$ is
complete. The following result is called Newman's Lemma
\cite{TAT}.

\begin{lem}
A Noetherian reduction relation on a set is confluent if it is
locally confluent.
\end{lem}

Hence, if $R$ is compatible with an admissible well-ordering on
$T$ and  $\to_R$ is locally confluent then $\to_R$ is complete. By
orienting the pairs of $R$ with respect to the chosen ordering
$>_T$ on $T$, $R$ is made to be Noetherian. The problem remaining
is testing for local confluence of $\to_R$ and changing $R$ in
order to obtain an equivalent confluent reduction relation.\\

We will now explain the notion of critical pair for a rewrite
system for $T$, extending the traditional notion to our situation.
In particular the overlaps involve either just $R_T$, or just
$R_P$ or an interaction between $R_T$ and $R_P$.

\begin{Def}
A term $crit\in T$ is called \textbf{critical} if it may be reduced by
two or more different rules
i.e. $crit \to_R crit1$, $crit \to_R crit2$ and $crit1 \not= crit2$.
A pair $(crit1,crit2)$ of distinct terms resulting from two single-step
reductions of the same term is called a \textbf{critical pair}.
A critical pair for a reduction relation $\to_R$ is said to
\textbf{resolve} if
there exists a (common) term $res$ such that both $crit1$ and $crit2$ reduce to a
$res$, i.e. $crit1 \stackrel{*}{\to}_R res$,
$crit2 \stackrel{*}{\to}_R res$.
\end{Def}

We now define overlaps of rules for our type of rewrite system, and show
how each kind results in a critical pair of the reduction relation.

If $t=x|b_1 \cdots b_n$, then a \textbf{part} of $t$ is either a term
$x|b_1 \cdots b_i$ for some $1 \leq  i \leq  n$ or a word $b_i b_{i+1} \cdots b_j$
for some $1 \leq  i \leq  j \leq  n$.

\begin{Def}
Let $(rule1,rule2)$ be a pair of rules of the rewrite
system $R=(R_T,R_P)$ where $R_T \subseteq T \times T$ and $R_P
\subseteq \arr\bP \times \arr \bP$.
If $rule1$ and $rule2$ may both be applied to the same term $crit$ in such a way
there is a part of the term $crit$ that
is affected by both the rules then we say that an \textbf{overlap} occurs.
\end{Def}

 There are five types of overlap for this kind of rewrite system,
 as shown in the following table:
 \vspace{2ex}
\begin{center}
\setlength{\extrarowheight}{6pt}
\begin{tabular}{||c|c|c|c|c|l|l||}
\hline
 \# &rule1 & in & rule2& in  & \hspace{1cm} overlap & critical pair \\
\hline
(i) &$(s_1,u_1)$ & $ R_T$  & $(s_2,u_2)$& $ R_T$  &  $s_2=s_1\cdot q$ for some $q \in \arr \bP$ & $(u_1\cdot q,u_2 )$ \\
\hline
(ii)&$(l_1, r_1)$ & $R_P$ & $(l_2,r_2)$ & $ R_P$   &  $l_1=pl_2q$ for some $p,q \in \arr \bP$    &  $(r_1,pr_2q)$
\\ \cline{1-1} \cline{6-7}
(iii) & & &  &   &  $l_1q=pl_2$  for some $p,q \in \arr \bP$   &   $(r_1q,pr_2)$
\\ \hline
(iv)& $(s_1,u_1)$ & $ R_T$ & $(l_1, r_1)$ & $R_P$  &  $s_1\cdot q= s\cdot l_1$ for some $s \in T, q \in \arr \bP$  & $(u_1\cdot q, s \cdot r_1)$
 \\\cline{1-1}  \cline{6-7}
 (v)&  & & &  &   $s_1 =s \cdot (l_1q)$ for some $s \in T, q \in \arr \bP$      & $(u_1, s \cdot r_1q)$  \\
  \hline
\end{tabular}

Overlap table

\end{center}
\vspace{2ex}

A pair of rules may overlap in more than one way, giving more than
one critical pair. For example the rules $(x|a^2ba, y|ba)$ and
$(a^2,b)$ overlap with critical term $x|a^2ba$ and critical pair
$(y|ba,x|b^2a)$ and also with critical term $x|a^2ba^2$ and
critical pair $(y|ba^2,x|a^2b^2)$.

\begin{lem}
Let $R$ be a finite rewrite system on the $\bP$-set $T$.
Consider applications of rules $rule1$ and $rule2$ affecting part $c$ of term
$t \in T$, resulting in a critical pair $(c_1,c_2)$ from $c$ and $(t_1,t_2)$
from $t$. If there is no overlap then $(t_1,t_2)$ resolves immediately.
Otherwise $(t_1,t_2)$ resolve providing $(c_1,c_2)$ does.
\end{lem}

\begin{proof}
Let $(t_1,t_2)$ be a critical pair. Then there exists a critical
term $t$ and two rules $rule1$, $rule2$ such that $t$ reduces to
$t_1$ with respect to $rule1$ and to $t_2$ with respect to
$rule2$.

We first give the two non-overlap cases.

Suppose $rule1:=(l_1,r_1)$, $rule2:=(l_2,r_2) \in R_P$.
Then there exist $s \in T$, $p, q \in \arr \bP$ such that
$t = s \cdot l_1 p l_2 q$ as shown:

$$
\xymatrix{\ar@{-}[rr]^*+{s} \ar@{ }[rr]_*+{s}
        & |
        & \ar@{-}@/^1pc/[r]^{r_1} \ar@{-}[r]_{l_1}
        & \ar@{-}[r]^p \ar@{ }[r]_p
        & \ar@{-}@/_1pc/[r]_{r_2} \ar@{-}[r]^{l_2}
        & \ar@{-}[r]^q \ar@{ }[r]_q
        &\\}
$$
The pair $(t_1,t_2)$ immediately resolves to $s \cdot r_1pr_2q$ by
applying $rule2$ to $t_1$ and $rule1$ to $t_2$.

Suppose that $rule1:=(s_1,u_1) \in R_T$ and $rule2:=(l_1,r_1)
\in R_P$ and the rules do not overlap.
Then there exist $p,q\in\arr\bP$ such that
$t = s_1 \cdot p l_1 q$ and then
$t_1 = u_1 \cdot p l_1 q$ and
$t_2 = s_1 \cdot p r_1 q$ as shown:
$$
\xymatrix{\ar@{-}[rr]_*+{s_1}
          \ar@{-}@/^1pc/[rr]^{u_1}
        & |
        & \ar@{-}[r]^p \ar@{ }[r]_p
        & \ar@{-}@/_1pc/[r]_{r_1} \ar@{-}[r]^{l_1}
        & \ar@{-}[r]^q \ar@{ }[r]_q
        &\\}
$$
The pair $(t_1,t_2)$ immediately resolves to $u_1 \cdot p r_1 q$ by
applying $rule2$ to $t_1$ and $rule1$ to $t_2$.

We now give the overlap cases in the order given in the table.

(i) Suppose $rule1:=(s_1,u_1), rule2:=(s_2,u_2) \in R_T$.
Then there exist $v,q\in\arr\bP$ such that
$c= s_1 \cdot q = s_2$, $t=c \cdot v$ and then $t_1 = u_1 \cdot q v$ and
$t_2 = u_2\cdot v$ as shown:
$$
\xymatrix{\ar@{-}[rr]
          \ar@{-}@/^1pc/[rr]^{u_1}
          \ar@{-}@/_1pc/[rrr]_{u_2}
        & | & \ar@{-}[r]^q
        & \ar@{-}[r]^v \ar@{ }[r]_v &\\}
$$
The critical pair here is
$(u_1 \cdot q, u_2)$ and if this resolves to $r$ then
$(t_1,t_2)$
resolves to $r \cdot v$.

Suppose $rule1:=(l_1,r_1)$, $rule2:=(l_2,r_2) \in R_P$.
There are  two possible overlap cases.

(ii) In the first case there exist
 $s \in T$, $p,q,v \in \arr \bP$ such that
$c = l_1 = p l_2 q$ and $t=s \cdot c v$ and then
$t_1 = s \cdot r_1 v$ and
$t_2 = s \cdot p r_2 q v$.
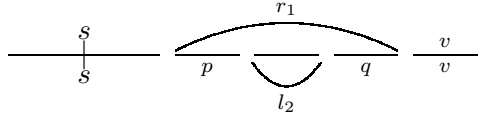
$$
\xymatrix{\ar@{-}[rr]^*+{s} \ar@{ }[rr]_*+{s}
        & | & \ar@{-}@/^1pc/[rrr]^{r_1} \ar@{-}[r]_p
        & \ar@{-}@/_1pc/[r]_{l_2} \ar@{-}[r]
        & \ar@{-}[r]_q
        & \ar@{-}[r]^{v} \ar@{ }[r]_v
        &\\}
$$
The critical pair here is
$(r_1, p r_2 q)$ and if this resolves to $r$ then
$(t_1,t_2)$ resolves to $s \cdot r v$.

(iii) In the second case there exist
$s \in T$, $p, q, v \in \arr \bP$ such that
$c = l_1 q = p l_2$ and $t= s \cdot cv$ and then
$t_1 = s \cdot r_1 q v$ and
$t_2 = s \cdot p r_2 v$.
$$
\xymatrix{\ar@{-}[rr]^*+{s}  \ar@{ }[rr]_*+{s}
        & |
        & \ar@{-}@/^1pc/[rr]^{r_1} \ar@{-}[r]_p
        & \ar@{-}@/_1pc/[rr]_{r_2} \ar@{-}[r]
        & \ar@{-}[r]^q
        & \ar@{-}[r]^{v} \ar@{ }[r]_v
        &\\}
$$
The critical pair is
$(r_1 q,p r_2)$ and if this resolves to $r$ then
$(t_1,t_2)$
resolves to $s\cdot r v$.

Suppose finally that $rule1:=(s_1,u_1) \in R_T$ and $rule2:=(l_1,r_1)
\in R_P$. Then there are  two possible overlap cases.

(iv) In the first case there exist
$s \in T$, $q, v \in \arr \bP$ such that
$c= s_1  = s \cdot l_1 q $ and $t=c \cdot v$ and then
$t_1 = u_1 v$ and
$t_2 = s r_1 q v$.
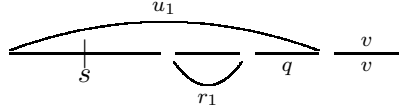
$$
\xymatrix{\ar@{-}[rr]_*+{s}
          \ar@{-}@/^1pc/[rrrr]^{u_1}
        & |
        & \ar@{-}@/_1pc/[r]_{r_1} \ar@{-}[r]
        & \ar@{-}[r]_q
        & \ar@{-}[r]^v \ar@{ }[r]_v &\\}
$$
The critical pair is
$(u_1, s \cdot r_1 q)$ and if this resolves to $r$
then $(t_1,t_2)$ resolves to $r \cdot v$.

(v) In the second case there exist
$s \in T$, $q, v \in \arr \bP$ such that
$c = s_1 \cdot q = s \cdot l_1 $  and $t=c \cdot v$ and then
$t_1 = u_1 \cdot q v$ and
$t_2 = s \cdot r_1 v$.
$$
\xymatrix{\ar@{-}[rr]_*+{s}
          \ar@{-}@/^1pc/[rrr]^{u_1}
        & |
        & \ar@{-}@/_1pc/[rr]_{r_1} \ar@{-}[r]
        & \ar@{-}[r]^q
        & \ar@{-}[r]^v \ar@{ }[r]_v
        &\\}
$$
The critical pair is
$(s_1 \cdot q, s \cdot r_1)$ and if this resolves to $r$
then $(t_1,t_2)$  resolves to $r \cdot v$.\\

Thus we have considered all possible ways in which a term may be
reduced by two different rules, and shown that resolution of the
critical pair (when not immediate) depends upon the resolution of
the critical pair resulting from a particular overlap of the
rules.
\end{proof}

\begin{cor} If all the overlaps between rules of a rewrite system $R$
on $T$ resolve then all the critical pairs for the reduction relation
$\to_R$ resolve, and so $\to_R$ is confluent.
\end{cor}
\begin{proof} This is immediate from the Lemma.
\end{proof}

\begin{lem}
All overlaps of a pair of rules of $R$ can be found by
looking for two types of overlap between the lists representing the left
hand sides of rules.
\end{lem}
\begin{proof}
Let $rule1=(l_1,r_1)$ and $rule2=(l_2,r_2)$ be a pair of rules. Recall
that
$\mathtt{List}(t)$
is the representation of a term $t \in T$ as a list.
The first type of list overlap occurs when $\mathtt{List}(l_2)$ is a
sublist of $\mathtt{List}(l_1)$ (or vice-versa). This happens in cases
(i), (ii) and (v).
The second type of list overlap occurs when the end of
$\mathtt{List}(l_1)$
matches
the beginning of $\mathtt{List}(l_2)$ (or vice-versa). This happens in
cases (iii) and (iv).
\end{proof}

The program for finding overlaps and the resulting critical pairs is
outlined in the algorithm below.

\begin{alg}(Critical Pairs)
Given a rewrite system $R$ all critical pairs are determined.
\begin{enumerate}
\item
(Input:) A rewrite system $R$ as a set of rules (pairs of lists).
\item
(Initialise:) Set $CRIT:=\emptyset$.
\item
(Procedure:) Take pairs of rules $(l_1,r_1)$ and $(l_2,r_2)$from $R$.
 Test (a) whether $\mathtt{List}(l_2)$ is a sublist of $\mathtt{List}(l_1)$.
If it is then find $u$ and $v$ such that $u \cdot l_2 v = l_1$.
Add the critical pair $(u \cdot r_2 v, r_1 )$ to $CRIT$.
       Now test (b) whether for $i=1,2 \ldots$ the sublist of length $i$ at the
       right of $\mathtt{List}(l_1)$ is equal to the sublist of length $i$ on the
       left of $\mathtt{List}(l_2)$.
       For each $i$ where this occurs, set $u$ to be the
       part of $\mathtt{List}(l_1)$ not in the overlap, and $v$ to be the part
       of $\mathtt{List}(l_2)$ not
       in the overlap.
       Add the critical pair
$( r_1 \cdot v, u \cdot r_2)$ to $CRIT$.
       Repeat the procedure until all (ordered) pairs of rules have been examined for
       overlaps.
\item
(Output:) An exhaustive list of critical pairs $CRIT$.
\end{enumerate}
\end{alg}

It has now been proved that all the critical pairs of a finite rewrite
system $R$ on $T$ can be listed.
To test whether a critical pair resolves, each side of it is reduced
using the function $\mathtt{Reduce}$. If $\mathtt{Reduce}$ returns the
same term for each side then the pair resolves.

\subsection{Completion Procedure}

We have shown: (i) how to find overlaps between rules of $R$;
(ii) how to test whether
the resulting critical pairs resolve; and
(iii) that if all the critical pairs resolve
then this imples $\to_R$ is confluent. We now show that critical pairs
which do not resolve may be added to $R$ without affecting the
equivalence relation $R$ defines on $T$.

\begin{lem}
Any critical pair $(t_1,t_2)$ of a rewrite system $R$
may be added to the rewrite system without
changing the equivalence relation ${\stackrel{*}{\lra}}_R$.
\end{lem}
\begin{proof}
By
definition $(t_1,t_2)$ is the result of two different single-step
reductions being
applied to a critical term $t$. Therefore $t\to_Rt_1$ and $t\to_Rt_2$.
It is
immediate that
$t_1 \, {\stackrel{*}{\lra}}_R \, t \, {\stackrel{*}{\lra}}_R \, t_2$,
and so
adding $(t_1,t_2)$ to $R$
does not add anything to the equivalence relation
${\stackrel{*}{\lra}}_R$.
\end{proof}

We have now set up and proved everything necessary for a variant of the
Knuth-Bendix procedure, which will add rules to a rewrite system $R$
resulting from a presentation of a Kan extension, to attempt to find an
equivalent complete rewrite system $R^C$. The benefit of such a system is that
$\to_{R^C}$ then acts as a normal form function for $\stackrel{*}{\lra}_{R^C}$
on $T$.

\begin{thm}
\label{proc} Let $\mathcal{P}= \lan \Gamma | \Delta | RelB |X |F
\ran$ be a finite presentation of a Kan extension $(K, \ep)$. Let
$P:=P \Delta$, $T:= \bigsqcup_{\ob\Delta} \bigsqcup_{\ob\Gamma} XA \times
\bP(FA, B)$, and let $R$ be the initial rewrite system for
$\mathcal{P}$ on $T$. Let $>_T$ be an admissible well-ordering on
$T$. Then there exists a procedure which, if it terminates,  will
return a rewrite system $R^C$ which is complete with respect to
the ordering $>_T$ and  such that the equivalence relations
$\stackrel{*}{\lra}_R$, $\stackrel{*}{\lra}_{R^C}$ coincide.
\end{thm}

\begin{proof}
The procedure finds all critical pairs resulting from overlaps of
rules of $R$. It attempts to resolve them. When they do not
resolve it adds them to the system as new rules. Critical pairs of
the new system are then examined. When all the critical pairs of a
system resolve, then the procedure terminates, the final rewrite
system $R^C$ obtained is complete. This procedure has been
verified in the preceding results of this section.
\end{proof}

\begin{alg}[Completion]

Given the presentation of a
Kan extension and the ordering $>_T$, a
complete rewrite system with respect to $>_T$
is determined -- if the algorithm terminates.
\begin{enumerate}
\item
(Input:) A rewrite system $R$ on $T$ and an ordering $>_T$ on $T$.
\item
(Initialise:) Set $NewRules:=R$ and $OldRules:=\emptyset$.
\item
(Loop:) While $NewRules \not= OldRules$, set $OldRules:=NewRules$.
   Use the algorithm {\sf Critical Pairs} to
   determine all the critical pairs of $NewRules$.
   Remove each critical pair in turn from the list,
   and reduce both sides of the pair with respect to $NewRules$ using
   the algorithm {\sf Reduce}.
   If the left entry is greater than the right (with respect
   to $>_T$) then add the reduced critical pair to $NewRules$.
   If the right entry is greater than the left then
   add the reversed, reduced critical pair to $NewRules$.
   Repeat this loop until all critical pairs resolve and no rules are added.
\item
(Output:) A complete rewrite system $NewRules$ on $T$.
\end{enumerate}
\end{alg}

Supposing that the completion procedure outlined above terminates, we
will now briefly discuss how to interpret the complete
rewrite system on $T$.

\section{Interpreting the Output}

\subsection{Finite Enumeration of the Kan Extension}

When every set $KB$ is finite we may catalogue the elements of all
of the sets $\bigsqcup_B  KB$ in stages.

The first stage catalogues  the elements $x|\id_{ FA}$ where $x
\in XA$ for some $A \in \ob \Gamma$. These elements are considered
to have length one. The next stage builds on the set of
irreducible elements from the last block to construct elements of
the form $x|b$ where $b:FA\to B$ for some $B \in \ob \Delta$. This
is effectively acting on the sets with the generating arrows to
define new (irreducible) elements of length two. The next stage
builds on the irreducibles from the last block by acting with the
generators again. When all the elements of a block of elements of
the same length are reducible then the enumeration terminates (any
longer term will contain one of these terms and therefore be
reducible). The set of irreducibles is a set of normal forms for
$\bigsqcup_B  KB$. The subsets $KB$ of $\bigsqcup_B  KB$ are
determined by the function $\bar{\tau}$, i.e. if $x|b_1 \cdots
b_n$ is a normal form in $\bigsqcup_B  KB$ and $\tau(x|b_1 \cdots
b_n):=tgt(b_n)=B_n$ then $x|b_1 \cdots b_n$ is a normal form in
$KB_n$. Of course if one of the sets $KB$ is infinite then this
may prevent the enumeration of other finite sets $KB_i$. The same
problem would obviously prevent a Todd-Coxeter completion. This
cataloguing method only applies to finite Kan extensions. It has
been implemented in the function $kan$, which has an enumeration
limit of 1000 set in the
program.

\subsection{Regular Expression for the Kan Extension}

Let $R$ be a finite complete rewrite system on $T$ for the Kan
extension $(K,\ep)$. Then the theory of languages and regular
expressions may be applied. The set of irreducibles in $T$ is
found after the construction of an automaton from the rewrite
system and the derivation of a language from this automaton.
Details of this method may be found in chapter four of \cite{Anne}.

\subsection{Iterated Kan Extensions}
One of the pleasant features of this procedure is that the input
and output are of similar form. The consequence of this is that if
the extended action $K$ has been defined on $\Delta$ then given a
second functor $G':\bB\to\bC$ and a presentation $cat\lan
\Lambda|RelC \ran$ for $\bC$ it is straightforward to consider a
presentation for the Kan extension data $(K',G')$. This new
extension is in fact the Kan extension with data $(X',G'\circ F')$

\begin{lem}
Let $kan\lan \Gamma|\Delta|RelB|X|F\ran$ be a presentation for a
Kan extension $(K,\ep)$. Let $cat\lan\Lambda|RelC\ran$ present a
category $\bC$ and let $G':\bB\to\bC$ be a functor. Then the Kan
extension presented by $kan\lan \Gamma|\Lambda|RelC|X|G\circ
F|\ran$ is equal to the Kan extension presented by $kan\lan
\Delta|\Lambda|RelC|K|G\ran$.
\end{lem}

\begin{proof}
Let $kan\lan \Gamma|\Delta|RelB|X|F\ran$ present the Kan extension
data $(X',F')$ for the Kan extension $(K,\ep)$. Let $\bC$ be a
category finitely presented by $cat\lan\Lambda|RelC\ran$ and let
$G':\bB\to\bC$. Then $kan\lan \Delta|\Lambda|RelC|K|G\ran$
presents the Kan extension data $(K',G')$  for the Kan extension
$(L,\eta)$.

We require to prove that $(L,\eta\circ \ep)$ is the Kan extension
presented by $kan\lan \Gamma|\Lambda|RelC|X|G\circ F\ran$ having
data $(X',G'\circ F')$. It is clear that $(L,\eta \circ \ep)$
defines an extension of the action $X$ along $G\circ F$ because
$L$ defines an action of $\bC$ and $\eta \circ \ep :X\to L\circ
G\circ F$ is a natural transformation. \\ For the universal
property, let $(M,\nu)$ be another extension of the action $X$
along $F\circ G$. Then consider the pair $(M\circ G,\nu)$, it is
an extension of $X$ along $F$. Therefore there exists a unique
natural transformation $\alpha:X \to M\circ G\circ F$ such that
$\alpha \circ \ep=\nu$ by universality of $(K,\ep)$. Now consider
the pair $(M,\alpha)$, it is an extension of $K$ along $G$.
Therefore there exists a unique natural transformation $\beta:L\to
M$ such that $\beta\circ\eta=\alpha$ by universality of
$(L,\eta)$. Therefore $\beta$ is the unique natural transformation
such that $\beta \circ \eta \circ \ep=\nu$, which proves the
universality of the extension $(L,\eta\circ\ep)$.
\end{proof}

\section{Example of a {\sf GAP} session on the Rewriting Procedure}
\label{kaneg}

Here we give an example to show the use of the implementation.
Let $\bA$ and $\bB$ be the categories generated by the graphs
below, where $\bB$ has the relation $b_1b_2b_3=b_4$. {\large
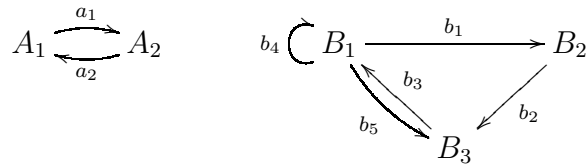
$$
\xymatrix{A_1 \ar@/^/[r]^{a_1}
        & A_2 \ar@/^/[l]^{a_2}
       && B_1  \ar@(dl,ul)^{b_4}
             \ar[rr]^{b_1}
             \ar@/_/[dr]_{b_5}
       && B_2 \ar[dl]^{b_2} \\
     &&&& B_3 \ar[ul]_{b_3} \\}
$$}
Let $X:\bA \to \sets$ be defined by $ XA_1 = \{ x_1, x_2, x_3 \},
                                   \ XA_2 = \{ y_1, y_2 \}$   with\\
$Xa_1:XA_1 \to XA_2: x_1 \mapsto y_1, x_2 \mapsto y_2, x_3 \mapsto
y_1$,\\ $Xa_2:XA_1 \to XA_2: y_1 \mapsto x_1, y_2 \mapsto x_2,$\\
and let $F:\bA \to \bB$ be defined by $FA_1=B_1, \ FA_2=B_2, \
Fa_1=b_1$ and $Fa_2 = b_3 b_2$. The input to the computer program
takes the following form.
First read in the program and set up the variables:
\begin{verbatim}
gap> RequirePackage("kan");
gap> F:=FreeGroup("b1","b2","b3","b4","b5","x1","x2","x3","y1","y2");;
gap> b1:=F.1;;b2:=F.2;;b3:=F.3;;b4:=F.4;;b5:=F.5;;
gap> x1:=F.6;;x2:=F.7;;x3:=F.8;;y1:=F.9;;y2:=F.10;;
\end{verbatim}
Then we input the data (choice of names is unimportant):
\begin{verbatim}
gap> OBJa:=[1,2];;
gap> ARRa:=[[1,2],[2,1]];;
gap> OBJb:=[1,2,3];;
gap> ARRb:=[[b1,1,2],[b2,2,3],[b3,3,1],[b4,1,1],[b5,1,3]];;
gap> RELb:=[[b1*b2*b3,b4]];;
gap> fOBa:=[1,2];;
gap> fARRa:=[b1,b2*b3];;
gap> xOBa:=[[x1,x2,x3],[y1,y2]];;
gap> xARRa:=[[y1,y2,y1],[x1,x2]];;
\end{verbatim}
To combine all this data in one record (the field names are important):
\begin{verbatim}
gap> KAN:=rec( ObA:=OBJa, ArrA:=ARRa,  ObB:=OBJb, ArrB:=ARRb, RelB:=RELb,
              FObA:=fOBa, FArrA:=fARRa, XObA:=xOBa, XArrA:=xARRa );;
\end{verbatim}
To calculate the initial rules do:
\begin{verbatim}
gap> InitialRules( KAN );
\end{verbatim}
The output will be:
\begin{verbatim}
i= 1, XA= [ x1, x2, x3 ], Ax= x1, rule= [ x1*b1, y1 ]
i= 1, XA= [ x1, x2, x3 ], Ax= x2, rule= [ x2*b1, y2 ]
i= 1, XA= [ x1, x2, x3 ], Ax= x3, rule= [ x3*b1, y1 ]
i= 2, XA= [ y1, y2 ], Ax= y1, rule= [ y1*b2*b3, x1 ]
i= 2, XA= [ y1, y2 ], Ax= y2, rule= [ y2*b2*b3, x2 ]
[ [ b1*b2*b3, b4 ], [ x1*b1, y1 ], [ x2*b1, y2 ], [ x3*b1, y1 ],
  [ y1*b2*b3, x1 ], [ y2*b2*b3, x2 ] ]
\end{verbatim}
This means that there are five initial $\ep$-rules:\\
$( \ x_1|Fa_1,  x_1.a_1|\id_{FA_2} \ ), \
 ( \ x_2|Fa_1,  x_2.a_1|\id_{FA_2} \ ),
( \ x_3|Fa_1,  x_3.a_1|\id_{FA_2} \ ),$\\ $
 ( \ y_1|Fa_2,  y_1.a_1|\id_{FA_1} \ ), \
 ( \ y_2|Fa_2,  y_2.|a_1 \id_{FA_1} \ ),$ \\
i.e. $ \ x_1|b_1 \to y_1|\id_{B_2}, \ x_2|b_1 \to y_2|\id_{B_2}, \
x_3|b_1 \to y_1|\id_{B_2}, \  y_1|b_2b_3  \to  x_1|\id_{B_1},
y_2|b_2b_3 \to x_2|\id_{B_1} \ $

and one initial $K$-rule: $b_1b_2b_3 \to  b_4$.

To attempt to complete the Kan extension presentation do:
\begin{verbatim}
gap> KB( KAN );
\end{verbatim}
The output is:
\begin{verbatim}
[ [ x1*b1, y1 ], [ x1*b4, x1 ], [ x2*b1, y2 ], [ x2*b4, x2 ], [ x3*b1, y1 ],
  [ x3*b4, x1 ], [ b1*b2*b3, b4 ], [ y1*b2*b3, x1 ], [ y2*b2*b3, x2 ] ]
\end{verbatim}
In other words to complete the system we have to add the rules
$$x_1|b_4  \to  x_1, \quad x_2|b_4  \to  x_2, \text{ and } x_3|b_4  \to
x_1.$$
The result of attempting to compute the sets by doing:
\begin{verbatim}
gap> Kan(KAN);
\end{verbatim}
is a long list and then:
\begin{verbatim}
enumeration limit exceeded: complete rewrite system is:
[ [ x1*b1, y1 ], [ x1*b4, x1 ], [ x2*b1, y2 ], [ x2*b4, x2 ], [ x3*b1, y1 ],
  [ x3*b4, x1 ], [ b1*b2*b3, b4 ], [ y1*b2*b3, x1 ], [ y2*b2*b3, x2 ] ]
\end{verbatim}
This means that the sets $KB$ for $B$ in $\bB$ are too large. The
limit set in the program is  1000. (To change this the user should type
{\tt EnumerationLimit:= 5000} -- or whatever, after reading in the program.)
In fact the above example is
infinite.
The complete rewrite system is output instead of the
sets. We can in fact use this to obtain regular expressions for
the sets. In this case the regular expressions are:
\begin{center}
\begin{tabular}{lcl}
$KB_1$ & $:=$ & $(x_1+x_2+x_3)|(b_5(b_3{b_4}^*b_5)^*b_3{b_4}^*+\id_{B_1}).$\\
$KB_2$ & $:=$ & $(x_1+x_2+x_3)|b_5(b_3{b_4}^*b_5)^*b_3{b_4}^*(b_1) +
(y_1+y_2)|\id_{B_2}.$\\
$KB_3$ & $:=$ & $(x_1+x_2+x_3)|b_5(b_3{b_4}^*b_5)^*(b_3{b_4}^*b_1b_2
+\id_{B_3}) + (y_1 + y_2)|b_2.$\\
\end{tabular}
\end{center}
The actions of the arrows are defined by concatenation followed by
reduction.\\ For example $x_1|b_5b_3b_4b_4b_5$ is an element of
$KB_3$, so $b_3$ acts on it to give $x_1|b_5b_3b_4b_4b_5b_3$ which
is irreducible, and an element of $KB_1$.

The general method of  obtaining regular expressions for these
computations will be given in a separate paper (see Chapter 4 of
\cite{Anne}).

\section{Special Cases of the Kan Rewriting Procedure}

\subsection{Groups and Monoids}

ORIGINAL PROBLEM: Given a monoid presentation $mon\lan
\Sigma|Rel\ran$, find a set of normal forms for the monoid
presented.\\
KAN INPUT DATA: Let $\Gamma$ be the graph with one
object and no arrows. Let $X0$ be a one point set. Let $\bB$
be generated by the graph $\Delta$ with one object and arrows
labelled by $\Sigma$, it has relations $Rel\bB$ given by the
monoid relations. The functor $F$ maps the object of $\Gamma$ to
the object of $\Delta$.\\
KAN EXTENSION: The Kan extension
presented by $kan\lan \Gamma|\Delta|RelB|X|F\ran$ is such that
$K0$ is a set of normal forms for the elements of the
monoid, the arrows of $\bB$ (elements of $PX$) act on the right of
$\bB$ by right multiplication. The natural transformation $\ep$
makes sure that the identity of $\bB$ acts trivially and helps to
define the normal form function. The normal form function is $w
\mapsto \ep_0(1)\cdot(w):=Kw(\ep_0(1))$.\\

In this case the method of completion is the standard Knuth-Bendix
procedure used for many years for working with monoid
presentations of groups and monoids. This type of calculation is
well documented.

\subsection{Groupoids and Categories}

ORIGINAL PROBLEM: To specify a set of normal forms for the
elements of a groupoid or category given by a finite category
presentation $cat\lan \Lambda|Rel \ran$.\\
KAN INPUT DATA: Let
$\Gamma$ be the discrete graph with no arrows and object set equal
to $\ob\Lambda$. Let $XA$ be a distinct one object set for each
$A\in\ob\Gamma$. Let $\bB$ be the category generated by
$\Delta:=\Lambda$ with relations $Rel\bB:=Rel$. Let $F$ be defined
by the identity map on the objects.\\
KAN EXTENSION: Then the Kan
extension presented by $kan\lan \Gamma | \Delta | RelB | X | F
\ran$ is such that $KB$ is a set of normal forms for the arrows of
the category with target $B$, the arrows of $\bB$ (elements of
$P\Gamma$) act on the right of $\bB$ by right multiplication. The
natural transformation $\ep$ makes sure that the identities of
$\bB$ act trivially and helps to define the normal form function.
The normal form function is $w \mapsto
\ep_A\cdot(w):=Kw(\ep_A)$.\\

\begin{example}\emph{
Consider the group $S_3$ presented by $\lan x,y|x^3,y^2,xyxy\ran.$
The elements are \\
$\{ \id, x, y, x^2, xy, yx \}$.
The covering groupoid is generated by the Cayley graph.}

{\Large
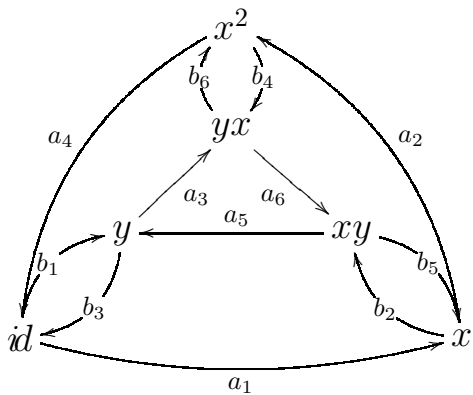
$$\xymatrix{&& x^2 \ar@/_1.7pc/[dddll]_{a_4} \ar@/^1pc/[d]|{b_4} && \\
            && yx \ar@/^1pc/[u]|{b_6} \ar[dr]_{a_6} && \\
            & y \ar@/^1pc/[dl]|{b_3} \ar[ur]_{a_3} &
            & xy \ar@/^1pc/[dr]|{b_5} \ar[ll]_{a_5} & \\
            \id  \ar@/^1pc/[ur]|{b_1} \ar@/_1pc/[rrrr]_{a_1} &&&&
            x \ar@/_1.7pc/[uuull]_{a_2} \ar@/^1pc/[ul]|{b_2} \\}
$$}

\emph{
The 12 generating arrows of the groupoid are $G \times X$:
$$
\{ [\id,x],[x,x],[y,x],\ldots,[yx,x],[\id,y],[x,y],\ldots,[yx,y] \}.
$$
To make calculations clearer, we relabel them
$\{ a_1, a_2, a_3, \ldots, a_6, b_1, b_2, \ldots,  b_6 \}$.
The groupoid has 18 relators (the boundaries of irreducible cycles of
the graph) $G \times R$,
the cycles may be written $[\id,x^3]$ and the corresponding boundary is
$[\id,x][x,x][x^2,x]$ i.e. $a_1a_2a_4$.
For the category presentation of the group we could add in the inverses
$\{ A_1, A_2, \ldots, A_6, B_1, B_2, \ldots, B_6 \}$ with the relators
$A_1a_1$
and $a_1A_1$ etc and end up with a category presentation with 24
generators and the 42 relations. In this case however the groupoid is
finite and so there is no need to do this. For example there would be no
need for $A_2$ because $(a_2)^{-1}=a_4a_1$.}

\emph{ Now suppose the left hand sides of two rules overlap (for
example $(a_1a_2a_4, \id)$ and $(a_4b_1a_3b_6, \id)$) in one of
the two possible ways previously described. Then we have a critical
pair $(b_1a_3b_6,a_1a_2)$ ). The following is {\sf GAP} output of the
completion of the rewrite system for the covering groupoid of our
example:}

\begin{verbatim}
gap> Rel;                                       ##Input rewrite system:
[ [ a1*a2*a4, IdWord ], [ a2*a4*a1, IdWord ], [ a4*a1*a2, IdWord ],
  [ a3*a6*a5, IdWord ], [ a6*a5*a3, IdWord ], [ a5*a3*a6, IdWord ],
  [ b1*b3, IdWord ], [ b3*b1, IdWord ], [ b2*b5, IdWord ],
  [ b5*b2, IdWord ], [ b4*b6, IdWord ], [ b6*b4, IdWord ],
  [ a1*b2*a5*b3, IdWord ], [ a2*b4*a6*b5, IdWord ],
  [ a3*b6*a4*b1, IdWord ], [ a4*b1*a3*b6, IdWord ],
  [ a5*b3*a1*b2, IdWord ], [ a6*b5*a2*b4, IdWord ] ]
gap> KB(Rel);                                   ##Completed rewrite
system:
[ [ b1*b3, IdWord ], [ b2*b5, IdWord ], [ b3*b1, IdWord ],
  [ b4*b6, IdWord ], [ b5*b2, IdWord ], [ b6*b4, IdWord ],
  [ a1*a2*a4, IdWord ], [ a1*a2*b4, b1*a3 ], [ a1*b2*a5, b1 ],
  [ a2*a4*a1, IdWord ], [ a2*a4*b1, b2*a5 ], [ a2*b4*a6, b2 ],
  [ a3*a6*a5, IdWord ], [ a3*a6*b5, b3*a1 ], [ a3*b6*a4, b3 ],
  [ a4*a1*a2, IdWord ], [ a4*a1*b2, b4*a6 ], [ a4*b1*a3, b4 ],
  [ a5*a3*a6, IdWord ], [ a5*a3*b6, b5*a2 ], [ a5*b3*a1, b5 ],
  [ a6*a5*a3, IdWord ], [ a6*a5*b3, b6*a4 ], [ a6*b5*a2, b6 ],
  [ b1*a3*a6, a1*b2 ],  [ b1*a3*b6, a1*a2 ], [ b2*a5*a3, a2*b4 ],
  [ b2*a5*b3, a2*a4 ],  [ b3*a1*a2, a3*b6 ], [ b3*a1*b2, a3*a6 ],
  [ b4*a6*a5, a4*b1 ],  [ b4*a6*b5, a4*a1 ], [ b5*a2*a4, a5*b3 ],
  [ b5*a2*b4, a5*a3 ],  [ b6*a4*a1, a6*b5 ], [ b6*a4*b1, a6*a5 ] ]
\end{verbatim}

\emph{
It is possible from this to enumerate elements of the category.
One method is to start with all the shortest arrows
($a_1,a_2,\ldots,b_6$) and see
which ones reduce and build inductively on the irreducible ones:\\
Firstly we have the six identity arrows
$\id_{\id}, \ \id_x, \ \id_y, \ \id_{x^2}, \ \id_{xy}, \ \id_{yx}$.\\
Then the generators $a_1, \ a_2, \ a_3, \ a_4, \ a_5, \ a_6, \ b_1, \
b_2, \
b_3, \ b_4, \ b_5, \ b_6$ are all irreducible.\\
Now consider paths of length 2:\\
$a_1a_2, \ a_1b_2, \ a_2a_4, \ a_2b_4, \
 a_3a_6, \ a_3b_6, \ a_4a_1, \ a_4b_1, \
 a_5a_3, \ a_5b_3, \ a_6a_5, \ a_6b_5, \
 b_1a_3, \ b_1b_3 \to \id_{\id},$\\
$b_2a_5, \ b_2b_5 \to \id_x, \ b_3a_1, \ b_3b_1 \to \id_y, \
 b_4a_6, \ b_4b_6 \to \id_{x^2}, \
 b_5a_2, \ b_5b_2 \to \id_{xy}, \ b_6a_4, \ b_6b_4 \to \id_{yx}$.\\
Building on the irreducible paths we get the paths of length 3:
$a_1a_2a_4 \to \id_{\id}, \ a_1a_2b_4 \to b_1a_3,$\\
$a_1b_2a_5 \to b_1, \ a_1b_2b_5 \to a_1, \ a_2a_4a_1 \to \id_x,\ldots$\\
All of them are reducible, and so we cannot build any longer paths; the
covering groupoid has 30 morphisms and 6 identity arrows and is the tree
groupoid with six objects.
}\end{example}

\begin{example}\emph{
This is a basic example to show how it is possible to specify the arrows
in an infinite small category with a finite complete presentation.
Let $\bC$ be the category generated by the following graph $\Gamma$
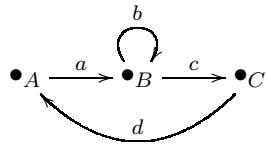
$$
\xymatrix{ \bullet_A \ar[r]^a & \bullet_B \ar@(ul,ur)^b \ar[r]^c
         & \bullet_C \ar@/^2pc/[ll]_d}
$$ with the relations $b^2c=c,\ ab^2=a$. This rewrite system is
complete, and so we can determine whether two arrows in the free
category $P\Gamma$ are equivalent in $\bC$. An automaton can be
drawn (see chapter 3 of \cite{Anne}), and from this we can
specify the language which is the set of normal forms. It is in
fact $$ a(cd(acd)*ab+bcd(acd)*ab) + b^{\dagger} + cd(acd)^*ab +
d(acd)^*ab $$ (and the three identity arrows) where $(acd)^*$ is
used to denote the set of elements of $\{ acd \}^*$ (similarly
$b^\dagger$), so $d(acd)*$, for example, denotes the set $\{d,
dacd, dacdacd, dacdacdacd, \ldots \}$, $+$ denotes the union and
$-$ the difference of sets. This is the standard notation for
languages and regular expressions. }\end{example}

\subsection{Coset systems and Congruences}

ORIGINAL PROBLEM:
Given a finitely presented group $G$ and a finitely generated subgroup
$H$
find a set of normal forms for the coset representatives
of $G$ with respect to $H$.\\
KAN INPUT DATA:
Let $\Gamma$ be the one object graph $\Gamma$ with arrows labelled by
the subgroup generators. Let $X0$ be a one point set on which the
arrows
of $\Gamma$ act trivially. Let $\bB$ be the category generated by the
one object
graph $\Delta$ with arrows labelled by the group generators, with the
relations
$Rel\bB$ of $\bB$ being the group relations. Let $F$ be defined on
$\Gamma$ by
inclusion of the subgroup elements to the group.\\
KAN EXTENSION:
The Kan extension presented by $kan\lan \Gamma|\Delta|RelB|X|F\ran$ is
such that
the
set
$K0$ is a set of representatives for the cosets, $Kb$ defines the
action
of the group on the cosets $Hg\mapsto Hgb$ and $\ep_0$ maps the
single
element of $X0$ to the representative for $H$ in $K0$.
Therefore it follows that the Kan extension defined is computable if and
only if
the coset system is computable. \\

In the monoidal case $F$ is the inclusion of the submonoid $\bA$ of the
monoid $\bB$, and the action is trivial as before. The Kan extension of
this
action gives the quotient of $\bB$ by the right congruence generated by
$\bA$,
namely the equivalence relation generated by $ab\sim
b$ for all $a \in \bA, b \in \bB$, with the induced right action of
$\bB$.\\

It is appropriate to give a calculated example here. The example is
infinite so
standard Todd-Coxeter methods will not terminate, but the Kan extension
/
rewriting procedures enable the complete specification of the coset
system.

\begin{example}\emph{
Let $\bB$ be the infinite group presented by
$$grp \lan a, b, c \ | \ a^2b=ba, a^2c=ca, c^3b=abc, caca=b \ran$$
and let $\bA$ be the subgroup generated by $\{ c^2 \}$. \\
We obtain one initial $\ep$-rule (because $\bA$ has one generating
arrow)
i.e. $ H|c^2 \to H|\id. $\\
We also have four initial $K$-rules corresponding to the relations for
$\bB$:}
$$a^2b \to ba, \ a^2c \to ca, \ c^3b \to abc, \ caca \to b.$$
\emph{
Note: On completion of this rewrite system for the group, we find 24
rules
and for all $n \in \bN$ both $a^n$ and $c^n$ are irreducibles with
respect to
this  system
(one way to prove the well-known fact that this the group is infinite).}\\

\emph{
The five rules are combined and an infinite complete system for the Kan
extension of the action is easily found (using Knuth-Bendix with the
length-lex order).
The following is the {\sf GAP} output of the set of 32 rules:}

\begin{verbatim}
[ [ H*b, H*a ], [ H*a^2, H*a ], [ H*a*b, H*a ], [ H*c*a, H*a*c ],
  [ H*c*b, H*a*c ], [ H*c^2, H ], [ a^2*b, b*a ], [ a^2*c, c*a ],
  [ a*b^2, b^2 ], [ a*b*c, c*b ], [ a*c*b, c*b ], [ b*a^2, b*a ],
  [ b*a*b, b^2 ], [ b*a*c, c*b ], [ b^2*a, b^2 ], [ b*c*a, c*b ],
  [ b*c*b, b^2*c ], [ c*a*b, c*b ], [ c*b*a, c*b ], [ c*b^2, b^2*c ],
  [ c*b*c, b^2 ], [ c^2*b, b^2 ], [ H*a*c*a, H*a*c ], [ H*a*c^2, H*a ],
  [ b^4, b^2 ], [ b^3*c, c*b ], [ b^2*c^2, b^3 ], [ b*c^2*a, b^2 ],
  [ c*a*c*a, b ], [ c^2*a^2, b*a ], [ c^3*a, c*b ], [ c*a*c^2*a, c*b ] ]
\end{verbatim}

\emph{
(Note that the rules without $H$ (i.e. the two-sided rules) constitute a
complete rewrite system for the group.)\\
The set $KB$ (recall that there is only one object $B$ of $\bB$)
is infinite. It is the set of (right) cosets of the subgroup in the
group.
Examples of these cosets include:
$$H, Ha, Hc, Ha^2, Hac, Ha^3, Ha^4, Ha^5,\ldots$$
A regular expression for the coset representatives is:
$$a^*+c+ac.$$
Alternatively consider the subgroup generated by $b$. Add the rule
$Hb \to H$ and the complete system below is obtained:}

\begin{verbatim}
[ [ H*a, H ], [ H*b, H ], [ H*c*a, H*c ], [ H*c*b, H*c ], [ H*c^2, H ],
  [ a^2*b, b*a ], [ a^2*c, c*a ], [ a*b^2, b^2 ], [ a*b*c, c*b ],
  [ a*c*b, c*b ], [ b*a^2, b*a ], [ b*a*b, b^2 ], [ b*a*c, c*b ],
  [ b^2*a, b^2 ], [ b*c*a, c*b ], [ b*c*b, b^2*c ], [ c*a*b, c*b ],
  [ c*b*a, c*b ], [ c*b^2, b^2*c ], [ c*b*c, b^2 ], [ c^2*b, b^2 ],
  [ b^4, b^2 ], [ b^3*c, c*b ], [ b^2*c^2, b^3 ], [ b*c^2*a, b^2 ],
  [ c*a*c*a, b ], [ c^2*a^2, b*a ], [ c^3*a, c*b ], [ c*a*c^2*a, c*b ] ]
\end{verbatim}
\emph{
(Again, the two-sided rules are the rewrite system for the group.)\\
This time the subgroup has index 2, and the coset representatives are
$\id$ and $c$.
}\end{example}

\subsection{Equivalence Relations and Equivariant Equivalence Relations}
ORIGINAL PROBLEM:
Given a set $\Omega$ and a relation $Rel$ on $\Omega$.
Find a set of representatives for the equivalence classes of the set
$\Omega$
under the equivalence relation generated by $Rel$.\\
KAN INPUT DATA:
Let $\Gamma$ be the graph with object set $\Omega$ and
generating arrows $a:A_1 \to A_2$ if $(A_1,A_2) \in Rel$.
Let $XA:=\{A\}$ for all $A\in\Omega$. The arrows of $\Gamma$ act
according to
the relation, so $src(a)\cdot a=tgt(a)$.
Let $\Delta$ be the graph with one object and no arrows so that $\bB$ is
the
trivial category with no relations. Let $F$ be the null functor.\\
KAN EXTENSION:
The Kan extension presented by $kan\lan \Gamma|\Delta|RelB|X|F\ran$
is such that $K0:= \Omega/{\stackrel{*}{\lra}}_{Rel}$ a set of
representatives for
the equivalence classes of the set $\Omega$ under the equivalence
relation generated by $Rel$.\\

Alternatively let $\Omega$ be a set with a group or monoid $M$ acting on
it.
Let $Rel$ be a relation on $\Omega$.
Define $\Gamma$ to have object set $\Omega$ and
generating arrows $a:A_1 \to A_2$ if $(A_1,A_2) \in Rel$ or if $A_1\cdot
m=A_2$
Again, $XA:=\{A\}$ for $A\in\ob\Gamma$ and the arrows act as in the case
above.
Let $\Delta$ be the one object graph with arrows labelled by generators
of $M$
and for $\bB$ let $Rel\bB$ be the set of monoid relations. Let $F$ be
the null
functor.
The Kan extension gives the action of $M$ on the
quotient of $X$ by the $M$-equivariant equivalence relation generated by
$Rel$.
This example illustrates the advantage of working in categories, since
this is a
coproduct of categories which is a fairly simple construction.

\subsection{Orbits of Actions}

ORIGINAL PROBLEM:
Given a group $G$ which acts on a set $\Omega$,  find a set $KB$ of
representatives for the orbits of the action of $\bA$ on $\Omega$.\\
KAN INPUT DATA:
Let $\Gamma$ be the one object graph with arrows labelled by the
generators of
the group. Let $X0:=\Omega$. Let $\Delta$ be the one object, zero
arrow
graph generating the trivial category $\bB$ with $Rel\bB$ empty. Let $F$
be the
null functor.\\
KAN EXTENSION:
The Kan extension presented by $kan\lan \Gamma|\Delta|RelB|X|F\ran$ is
such that
$K0$ is a set of representatives for the orbits of the
action of the group on $\Omega$.

We present a short example to demonstrate the procedure in this case.

\begin{example}\emph{
Let $\bA$ be the symmetric group on three letters with presentation\\
$mon\lan a,b| a^3,b^2,abab\ran$ and let $X$ be the set $\{v,w,x,y,z\}$.
Let
$\bA$ act on $X$ by giving $a$ the effect of the permutation
$(v \ w \ x)$ and $b$ the effect of $(v \ w)(y \ z)$.}\\

\emph{
In this calculation we have a number of $\ep$-rules and no $K$-rules.
The $\ep$-rules just list the action, namely (trivial actions omitted):
$$v\to w, w \to x, x \to v, v \to w, w \to v, y \to z, z \to y. $$
The system of rules is complete and reduces to $\{ w \to v, x \to v, z
\to y\}$.
Enumeration is simple: $v, \ w \to v, \ x \to v, \ y, \ z \to y$ so
there are
two orbits of $\Omega$ represented by $v$ and $y$.\\
This is a small example. With large
examples the idea of having a minimal element (normal form) in each
orbit
to act as an anchor or point of comparison makes a lot of sense.
This situation serves as another illustration of rewriting in the
framework of
a Kan extension, showing not only that rewriting gives a result, but
that it is
the procedure one uses naturally to do the calculation.
}\end{example}

One variation of this is if $\Omega$ is the set of elements of the group
and the
action is conjugation: $x^a := a^{-1}xa$. Then the orbits are the
\emph{conjugacy classes} of the group.

\begin{example}\emph{ Consider
the quarternion group, presented by $\lan a,b \ | \ a^4, b^4,
abab^{-1}, a^2b^2\ran$, and (we can enumerate the elements using the
variation
of
the Kan extensions method described in Example 3) $\Omega=\{ \id, a, b,
a^2, ab, ba, a^3, a^2b \}$. Construct the Kan extension as above, where
the actions of $a$ and $b$ are by conjugation on elements of $\bA$.\\
There are 16 $\ep$-rules which reduce to $\{ a^3 \to a, \ a^2b \to b, \
ba \to ab \}$. The conjugacy classes are enumerated by applying these
rules to the elements of $\bA$. The irreducibles are $\{\id, \ a, \ b, \
a^2, \ ab\}$, and these are representatives of the five conjugacy
classes.
}\end{example}

\subsection{Colimits of Diagrams of $\sets$}

ORIGINAL PROBLEM: Given a presentation of a category action
$act\lan \Gamma|X\ran$ find the colimit of the diagram in $\sets$
on which the category action is defined.\\ KAN INPUT DATA: Let
$\Gamma$ and $X$ be those given by the action presentation. Let
$\Delta$ be the graph with one object and no arrows that generates
the trivial category $\bB$ with $Rel\bB$ empty. Let $F$ be the null
functor.\\ KAN EXTENSION: The Kan extension presented by $kan\lan
\Gamma|\Delta|RelB|X|F\ran$ is such that $K0$ is the colimit
object, and $\ep$ is the set of colimit functions of the functor
$X:\bA \to \sets$.\\

Particular examples of this are when $\bA$ has two objects $A_1$
and $A_2$, and two non-identity arrows $a_1$ and $a_2$ from $A_1$
to $A_2$, \ and $Xa_1$ and $Xa_2$ are functions from the set $XA_1$
to the set $XA_2$ (\emph{coequaliser} of $a_1$ and $a_2$ in
$\sets$); $\bA$ has three objects $A_1$, $A_2$ and $A_3$ and two
non-identity arrows $a_1:A_1 \to A_2$ and $a_2:A_1 \to A_3$. \
$XA_1$, $XA_2$ and $XA_2$ are sets, and $Xa_1$ and $Xa_2$ are
functions between these sets (\emph{pushout} of $a_1$ and $a_2$ in
$\sets$). The following example is included not as an illustration
of rewriting but to show another situation where presentations of
Kan extensions can be used to express a problem naturally.

\begin{example}\emph{
Suppose we have two sets $\{ x_1, x_2, x_3 \}$ and $\{y_1, y_2,
y_3, y_4 \}$,
with two functions from the first to the second given by
$( x_1 \mapsto y_1,  x_2 \mapsto y_2, x_3 \mapsto y_3 )$ and
$( x_1 \mapsto y_1,  x_2 \mapsto y_1, x_3 \mapsto y_3 )$.\\
Then we can calculate the coequaliser.
We have a number of $\ep$-rules
$$
y_1|\id_0 \to x_1|\id_0,  y_2|\id_0 \to
x_2|\id_0,
y_3|\id_0 \to x_3|\id_0, y_1|\id_0 \to
x_1|\id_0,
y_2|\id_0 \to x_1|\id_0, y_3|\id_0 \to
x_3|\id_0.
$$
There is just one overlap, between $(y_2|\id_0 \to
x_1|\id_0)$ and
$(y_2|\id_0 \to x_2|\id_0)$: to
resolve the critical pair we add the rule $x_2|\id_0 \to
x_1|\id_0$,
and the system is complete:
$$\{y_1|\id_0 \to x_1\id_0|, \ y_2|\id_0 \to
x_1|\id_0,
\ y_3|\id_0 \to x_3|\id_0, \ x_2|\id_0 \to
x_1|\id_0
\}.$$
The elements of the set $K0$ are easily enumerated:\\
$$
x_1|\id_0, \ x_2|\id_0 \to x_1|\id_0, \
x_3|\id_0, \ y_1|\id_0 \to x_1|\id_0, \
y_2|\id_0 \to x_1|\id_0, \
y_3|\id_0 \to x_3|\id_0, \
y_4|\id_0.
$$
So the coequalising set is
$$K0=\{x_1|\id_0, x_3|\id_0, y_4|\id_0 \},$$
and the coequaliser function
to it from $XA_2$ is given by $y_i\mapsto y_i|\id_0$ for
$i=1,\ldots,4$
followed by reduction defined by $\to$ to an element of $K0$.
}\end{example}

\subsection{Induced Permutation Representations}
\label{appropact}

Let $\bA$ and $\bB$ be groups and let $F:\bA \to \bB $ be a
morphism of groups. Let $\bA$ act on the set $XA$. The Kan
extension of this action along $F$ is known as the action of $\bB$
{\em induced} from that of $\bA$ by $F$, and is written $F_*(XA)$.
It can be constructed simply as the set $X \times \bB$ factored by
the equivalence relation generated by $(xa,b)\sim (x,F(a)b)$ for
all $x \in XA, a \in \bA,b \in \bB$. The natural transformation
$\ep$ is given by $ x \mapsto [x,1]$, where $[x,b]$ denotes the
equivalence class of $(x,b)$ under the equivalence relation
$\sim$. The morphism $F$ can be factored as an epimorphism
followed by a monomorphism, and there are other descriptions of
$F_*(XA)$ in these cases, as follows.

Suppose first that $F$ is an epimorphism with kernel $N$. Then we can
take as a representative of $F_*(XA)$ the orbit set $X/N$ with the
induced action of $\bB$.

Suppose next that $F$ is a monomorphism, which we suppose is an
inclusion. Choose a set $T$ of representatives of the right cosets
of $\bA$ in $\bB$, so that $1 \in T$. Then the induced
representation can be taken to be $XA \times T$ with $\ep$ given
by $x \mapsto (x,1)$ and the action given by $(x,t)^b= (xa,u)$
where $t,u \in T, b \in \bB, a \in \bA$ and $tb=au$.

On the other hand, in practical cases, this factorisation of $F$
may not be a convenient way of determining the induced
representation.

In the case $\bA,\bB$ are monoids, so that $X$ is a transformation
representation of $\bA$ on the set $XA$, we have in general no
convenient description of the induced transformation representation
except by one form or another of the construction of the Kan
extension.
This yields a quotient of the free product of the monoids
$\{x\}\times\bB$, $x \in XA$ by the equivalence relation generated by
$(x,F(a)b) \sim (x\cdot a, b)$, $a \in \bA, b \in \bB$.

\end{document}